%% file: main_arxiv.tex
\documentclass{article}

\usepackage[utf8]{inputenc} % allow utf-8 input
\usepackage[T1]{fontenc}    % use 8-bit T1 fonts
\usepackage{url}            % simple URL typesetting
\usepackage{booktabs}       % professional-quality tables
\usepackage{amsfonts}       % blackboard math symbols
\usepackage{nicefrac}       % compact symbols for 1/2, etc.
\usepackage{microtype}      % microtypography
\usepackage{graphicx}
\usepackage{doi}
\usepackage{orcidlink}
\usepackage{amsthm}
\usepackage[left=1.5cm, right=2cm]{geometry}

\usepackage{bm}
\usepackage{algorithm}% http://ctan.org/pkg/algorithm
\usepackage{algpseudocode}
\usepackage{color}
\usepackage{systeme}
\usepackage{comment}
\usepackage{scalerel}
\usepackage{graphicx}
\usepackage{grffile}
\usepackage{multirow}
\usepackage{pdflscape}
\usepackage{amsmath,amssymb}
\usepackage{cases}
\usepackage{bm}
\usepackage{booktabs}
\usepackage{tabularx}   % X column type (auto-width)

\usepackage[capitalise,noabbrev]{cleveref}
\usepackage{algorithm}
\usepackage{algorithmicx}
\usepackage{algpseudocode}
\usepackage{mathtools}
\usepackage{comment}
\usepackage{tikz}
\newlength{\figureheight}
\newlength{\figurewidth}
\usepackage{tikzsymbols}
\usepackage{pgfplots}
\pgfplotsset{compat=newest} %use latest pgf version
\usetikzlibrary{plotmarks}
\usetikzlibrary{arrows.meta}
\usetikzlibrary{positioning}
\usepgfplotslibrary{patchplots}
\usepackage{grffile}
\usepackage{multirow}
\usepackage{pdflscape}
 \usepackage{fontawesome}
\pgfplotsset{every axis/.append style={
                    label style={font=\scriptsize},
                    tick label style={font=\scriptsize},
                    legend style={font=\scriptsize}
                    }}
\pgfplotsset{compat=newest}
\pgfplotsset{plot coordinates/math parser=false}
\pgfplotsset{grid style={dotted,gray}}
\usepgfplotslibrary{groupplots}
\pgfdeclarelayer{background}
\pgfdeclarelayer{foreground}
\pgfsetlayers{background,main,foreground}
\usepackage{subcaption}
\usepackage{textcomp}
\usepackage{xcolor}
\usepackage{dsfont}
\usepackage{markdown}
\usepackage{accents}

\usepackage{hyperref}       % hyperlinks

\newlength{\dhatheight}

\newlength{\dtildeheight}
\newcommand{\doubletilde}[1]{%
    \settoheight{\dtildeheight}{\ensuremath{\tilde{#1}}}%
    \addtolength{\dtildeheight}{-0.30ex}%
    \tilde{\vphantom{\rule{1pt}{\dtildeheight}}%
    \smash{\tilde{#1}}}}

\newcommand\norm[1]{\left\lVert#1\right\rVert}

\newtheorem{thrm}{Theorem}[section]

\newtheorem{rmrk}[thrm]{Remark}

\definecolor{darkorange25512714}{RGB}{255,127,14}
\definecolor{forestgreen4416044}{RGB}{44,160,44}
\definecolor{steelblue31119180}{RGB}{31,119,180}
\definecolor{crimson2143940}{RGB}{214,39,40}

\makeatletter
\patchcmd{\ALG@step}{\addtocounter{ALG@line}{1}}{\refstepcounter{ALG@line}}{}{}
\newcommand{\ALG@lineautorefname}{Line}
\makeatother

\usepackage{todonotes}

\begin{document}

\title{A first-optimize-then-reduce framework for optimal control with the shifted proper orthogonal decomposition}

\author{
    Tobias Breiten\,\orcidlink{0000-0002-9815-4897}\thanks{\texttt{\{tobias.breiten@,burela@tnt.,pschulze@math.\}tu-berlin.de}}\,\,\thanks{Institute of Mathematics, Technische  Universit\"at Berlin, Germany},
    Shubhaditya Burela\,\orcidlink{0009-0003-4442-5297}\footnotemark[1]\,\,\footnotemark[2],
    Philipp Schulze\,\orcidlink{0000-0002-7299-4628}\footnotemark[1]\,\,\footnotemark[2]
}

\providecommand{\keywords}[1]
{
  \small	
  \textbf{Keywords:} #1\\
}
\providecommand{\ams}[1]
{
  \small	
  \textbf{AMS subject classifications:} #1
}

\maketitle

\begin{abstract}
Overcoming the Kolmogorov barrier for constructing efficient reduced-order models (ROMs) for transport-dominated problems remains a challenge. This has impeded their use as computationally cheap surrogates for partial differential equations (PDEs) associated with optimal control problems. Since such problems require multiple computations of the full PDE, employing their reduced-order surrogates instead could speed up the overall optimal control problem. Motivated by this idea, in this paper we explore the use of a nonlinear model-order reduction technique, namely, the shifted proper orthogonal decomposition (sPOD) in an optimal control context. In doing so, we explore the framework of first-optimize-then-reduce (FOTR) where an optimality system for the full PDE problem is constructed first, followed by approximating the optimality system with reduced-order models. We consider this framework for a linear quadratic optimal control problem constrained by a 1D linear advection equation and compare the computational performance of the sPOD method against the use of the standard POD.
\end{abstract}
\\\noindent\ams{35L02, 49M41, 49K20, 35Q35}

\section{Introduction}
We consider the following optimal control problem with a standard quadratic tracking-type cost functional,
\begin{equation}\label{eq:FOM_costFunc_continuous}
    \underset{u \in L^2(0, T; U)}{\mathrm{min}} \;\mathcal{J}(q, u) := \frac{1}{2} \int^{T}_0 \norm{q(t) - q_\mathrm{d}(t)}^2_H \mathrm{d}t + \frac{\mu}{2} \int_0^{T}\norm{u(t)}^2_{U}\,\mathrm{d}t\, ,
\end{equation} 
for $H:= L^2(0, l)$, $U:= \mathbb{R}^m$, $q_\mathrm{d} \in L^2(0, T; H)$, and $\mu>0$ such that
\begin{equation}\label{eq:PDE_continuous}
\begin{aligned}
    \tfrac{\mathrm{d}}{\mathrm{d} t} q(t) &= \mathcal{A} q(t) + \mathcal{B}u(t),  &  \quad t \in (0,T],\\
    q(0) &= q_0,
\end{aligned}
\end{equation}
where $\mathcal{A} = -v \frac{\mathrm{d}}{\mathrm{d}x}$, $\mathcal{A}: \mathcal{D}(\mathcal{A}) \subset L^2(0, l) \rightarrow L^2(0, l)$ for a constant velocity $v\in\mathbb{R}$ is a densely defined and closed linear operator that generates a strongly continuous (semi)group, with $\mathcal{D}(\mathcal{A}) = H^1_\mathrm{per}(0, l):= \{f \in H^1(0, l) | f(0) = f(l)\}$ see, e.g., \cite[Examples 2.6.12 \& 2.7.12, Thm.~3.8.6, Cor.~2.1.8, Prop.~2.3.1]{tucsnak_observation_2009}.
Given the control operator $\mathcal{B}\in\mathcal{L}(U,H)$ with $H:=L^2(0,l), U = \mathbb{R}^{m}$, which is taken directly from \cite{breiten2026optimalcontrolshiftedproper} where we study the first-reduce-then-optimize framework, the term $\mathcal{B}u(t) := \sum^{m}_{k=1} b_k u_k(t)$ is the control term, $u_k(t) \in U$ being the control intensities and $b_k \in \mathcal{D}(\mathcal{A})$ the control shape functions.

To solve the above mentioned optimal control problem, we follow the formal Lagrange method \cite{troltzsch_optimal_2010} and accordingly formulate the Lagrangian,
\begin{equation*}
    \mathcal{L}(q, u, \lambda):= \mathcal{J}(q, u) - \int^T_0 \bigl\langle \lambda(t), \tfrac{\mathrm{d}}{\mathrm{d} t} q(t) - \mathcal{A} q(t) - \mathcal{B}u(t) \bigr\rangle_H \mathrm{d} t\,.
\end{equation*}
Subsequently, if $(\bar{q}, \bar{u})$ is a local minimizer of $\mathcal{J}$ subject to \eqref{eq:PDE_continuous}, then there exists $\lambda$ (the adjoint) such that the triplet $(\bar{q}, \bar{u}, \lambda)$ satisfies the following first-order necessary optimality conditions
\begin{subnumcases}{\mathrm{OC}_{\scaleto{\mathrm{FOM}}{4pt}} := }
    \tfrac{\mathrm{d}}{\mathrm{d} t} q(t) = \mathcal{A} q(t) + \mathcal{B} u(t),\quad t \in (0,T]\label{eq:OC_FOM_1},  \\
    y(0) = y_0\label{eq:OC_FOM_2}, \\[1ex]
    -\tfrac{\mathrm{d}}{\mathrm{d} t}\lambda(t) = \mathcal{A}^* \lambda(t) + q(t) - q_\mathrm{d}(t),\quad t \in [0,T)\label{eq:OC_FOM_3}, \\
    \lambda(T) = 0\label{eq:OC_FOM_4},\\[1ex]
    \mu u(t) + \mathcal{B}^* \lambda(t)  = 0,\quad t \in [0,T], \label{eq:OC_FOM_5} 
\end{subnumcases} 
where \eqref{eq:OC_FOM_1}, \eqref{eq:OC_FOM_2} is referred to as the full-order model (FOM) state equation, \eqref{eq:OC_FOM_3}, \eqref{eq:OC_FOM_4} is the FOM adjoint equation, and \eqref{eq:OC_FOM_5} is the optimality condition for the control. For the particular case of the linear advection equation, we have $\mathcal{A}^* = - \mathcal{A} = v \frac{\mathrm{d}}{\mathrm{d}x}$ with $\mathcal{D}(\mathcal{A}^*) = H^1_\mathrm{per}(0, l)$ see \cite[Section 2.8]{tucsnak_observation_2009}.  
Additionally, theoretical results, such as the well-posedness of the state and the adjoint equation and the existence of an optimal control for an $L^2$ distributed control, follow straightforwardly from \cite{troltzsch_optimal_2010}. 

Optimal control problems become computationally challenging to solve, especially when higher-dimensional (hyperbolic) problems are considered. Therefore, our aim is to approximate $\mathrm{OC}_{\scaleto{\mathrm{FOM}}{4pt}}$ with appropriate ROMs that can potentially speed up the overall computation without unduly compromising the accuracy of the problem. Popular ROMs for such applications are often built on conventional linear methods such as proper orthogonal decomposition (POD). Owing to easier construction and simplified theoretical analysis, the application of POD-based ROMs in optimal control problem has been studied extensively, see, e.g., \cite{ravindran_reduced-order_2000, kunisch_control_1999, hinze_proper_2005, Afanasiev_wake_2001, grasle_pod_2017}.  Some of the popular methods that have gained significant theoretical and applicative traction are the trust region POD (TRPOD) \cite{arian_trust-region_nodate} or the optimality system POD (OSPOD) \cite{kunisch_proper_2008}. For a more detailed description of the plethora of techniques, see \cite{Benner2014}. However, for systems exhibiting a large Kolmogorov-$n$ width such as the advection equation with initial values featuring sharp wave fronts, linear MOR methods do not provide low-dimensional and accurate ROMs. In such scenarios, nonlinear methods are used to circumvent this problem, see \cite{HesPU26} for a detailed review of such methods. In this work, we specifically use the sPOD method \cite{reiss_shifted_2018, reiss_optimization-based_2021, krah_non-linear_nodate, schulze_energy-based}. In contrast to our previous work \cite{breiten2026optimalcontrolshiftedproper} on a first-reduce-then-optimize (FRTO) framework, here we take a different route and use the FOTR framework instead.

\section{Theoretical background}
In this section, we derive the reduced-order approximations for $\mathrm{OC}_{\scaleto{\mathrm{FOM}}{4pt}}$ considering both the POD-Galerkin (POD-G) and the sPOD-Galerkin (sPOD-G) methods. While the approximations with the POD-Galerkin method have already been well studied (see \cite{negri_reduced_2013} and references therein), the reduced-order approximations of $\mathrm{OC}_{\scaleto{\mathrm{FOM}}{4pt}}$ with the sPOD-Galerkin method are non-trivial and to the best of our knowledge have not been studied yet, except in a different framework in \cite{breiten2026optimalcontrolshiftedproper}.

\subsection{POD-Galerkin approximation of $\mathrm{OC}_{\scaleto{\mathrm{FOM}}{4pt}}$}
By solving the FOM state equation \eqref{eq:OC_FOM_1}, \eqref{eq:OC_FOM_2}, let us generate a snapshot set of the obtained trajectories $q \in C([0, T]; \mathcal{D}(\mathcal{A}))$ as
\begin{equation*}
    \mathcal{Q} := \text{span} \Bigl\{\, q(t_j) \ | \ t_j \in [0, T]   \text{ for } \: 1 \leq j \leq n_t \Bigr\} \subset \mathcal{D}(\mathcal{A})\;.
\end{equation*}
The POD attempts to approximate the trajectories as linear combination of orthonormal spatial basis functions $\{\phi_i\}^{\ell}_{i=1}$ (often called the POD basis or POD modes) and their respective time coefficients $\alpha_i(t) = \langle q(t), \phi_i\rangle_H $ (often called the reduced states) and results in
\begin{equation}\label{eq:Galerkin_POD_ansatz}
    q(t) \approx \sum^{\ell}_{i=1} \alpha_i(t)\phi_i,
\end{equation} 
by solving the following optimization problem,
\begin{subnumcases}{} 
    \mathrm{min} \: \frac{1}{2} \int^{T}_0 \norm{q(t) - \sum^{\ell}_{i=1} \Bigl\langle q(t), \phi_i \Bigr \rangle_H \: \phi_i}^2_H \mathrm{d}t ,\\[1em]
    \text{s.t.} \: \{\phi_i\}^{\ell}_{i=1} \subset \mathcal{D}(\mathcal{A}) \: \text{and} \:\langle \phi_i, \phi_j \rangle_H = \delta_{ij}, \: i, j = 1, \ldots , \ell \;.
\end{subnumcases}
Once the POD basis is constructed, one can derive the time-dynamical POD-G ROM by using the POD ansatz \eqref{eq:Galerkin_POD_ansatz} in \eqref{eq:OC_FOM_1},\eqref{eq:OC_FOM_2} and subsequently using the Galerkin orthogonality condition to obtain,
\begin{equation*}
\begin{aligned}
   &\Bigl \langle \phi_j, \sum^{\ell}_{i=1} \dot{\alpha}_i(t) \phi_i \Bigr \rangle_H = \Bigl \langle \phi_j, \sum^{\ell}_{i=1}\alpha_i(t)\mathcal{A}\phi_i \Bigr \rangle_H  +\Bigl\langle \phi_j, \mathcal{B}u(t) \Bigr \rangle_H \\
    &\Bigl \langle \phi_j, \sum^{\ell}_{i=1}\alpha(0)\phi_i \Bigr \rangle_H  = \Bigl\langle \phi_j, q_0 \Bigr\rangle_H
\end{aligned}
\end{equation*}
which can be compactly written as
\begin{equation}\label{eq:POD_Galerkin_continuous}
\begin{aligned}
   & \dot{\alpha}(t) = A_{\ell} \alpha(t) + B_{\ell} u(t)\\
    &\alpha(0)  = \alpha_0
\end{aligned}
\end{equation}
where $\quad \alpha_0 = [\langle \phi_j, q_0\rangle_H ]^{\ell}_{j=1} \in \mathbb{R}^{\ell}$ and 
\begin{equation*} 
\begin{aligned}
    B_{\ell} := \begin{bmatrix} \langle \phi_1, \mathcal{B}e_1\rangle_H  & \ldots  & \langle \phi_1, \mathcal{B}e_m\rangle_H \\ \vdots & \ddots & \vdots \\
     \langle \phi_{\ell}, \mathcal{B}e_1\rangle_H & \ldots  & \langle \phi_{\ell}, \mathcal{B}e_m\rangle_H \end{bmatrix}& \in \mathbb{R}^{\ell \times m},
     \quad A_\ell := \begin{bmatrix} \langle \phi_1, \mathcal{A}\phi_1\rangle_H  & \ldots  & \langle \phi_1, \mathcal{A}\phi_{\ell}\rangle_H \\ \vdots & \ddots & \vdots \\
     \langle \phi_{\ell}, \mathcal{A}\phi_1\rangle_H & \ldots  & \langle \phi_{\ell}, \mathcal{A}\phi_{\ell}\rangle_H \end{bmatrix} \in \mathbb{R}^{\ell \times \ell}\;. \\[1ex]
     \end{aligned}
\end{equation*}
For the adjoint equation \eqref{eq:OC_FOM_3}, \eqref{eq:OC_FOM_4}, we follow a similar derivation although with the POD ansatz for the adjoint variable $\lambda$,
\begin{equation}\label{eq:Galerkin_POD_ansatz_adjoint}
    \lambda(t) \approx \sum^{\ell_a}_{i=1} \lambda^{\ell_a}_i(t)\psi_i,
\end{equation} 
with orthonormal POD modes $\psi_1,\ldots,\psi_{\ell_a}$ for the adjoint state, and obtain
\begin{equation}
\begin{aligned}
    -\dot{\lambda}^{\ell_a}(t) &= A_{\ell_a} \lambda^{\ell_a}(t) + M \alpha(t) - \hat{q}_\mathrm{d}(t),  \\
    \lambda^{\ell_a}(T) &= 0,
\end{aligned}
\end{equation}
where
\begin{equation*}
    M := \begin{bmatrix} \langle \psi_1, \phi_1\rangle_H  & \ldots  & \langle \psi_1, \phi_{\ell}\rangle_H \\ \vdots & \ddots & \vdots \\
     \langle \psi_{\ell_a}, \phi_1\rangle_H & \ldots  & \langle \psi_{\ell_a}, \phi_{\ell}\rangle_H \end{bmatrix} \in \mathbb{R}^{\ell_a \times \ell}\;, \quad A_{\ell_a} := \begin{bmatrix} \langle \psi_1, \mathcal{A}^*\psi_1\rangle_H  & \ldots  & \langle \psi_1, \mathcal{A}^*\psi_{\ell_a}\rangle_H \\ \vdots & \ddots & \vdots \\
     \langle \psi_{\ell_a}, \mathcal{A}^*\psi_1\rangle_H & \ldots  & \langle \psi_{\ell_a}, \mathcal{A}^*\psi_{\ell_a}\rangle_H \end{bmatrix} \in \mathbb{R}^{\ell_a \times \ell_a}
\end{equation*}
with $\hat{q}_\mathrm{d}(t) = [\langle \psi_j, q_\mathrm{d}(t)\rangle_H ]^{\ell_a}_{j=1} \in \mathbb{R}^{\ell_a}$.
Lastly, the reduced-order approximation of the optimality condition for the control is given as
\begin{equation*}
    \mu u(t) + B_{\ell_a}^\top \lambda^{\ell_a}(t) = 0
\end{equation*}
where
\begin{equation*}
    B_{\ell_a} := \begin{bmatrix} \langle \psi_1,\mathcal{B}e_1 \rangle_H  & \ldots  & \langle \psi_{1}, \mathcal{B}e_m \rangle_H \\ \vdots & \ddots & \vdots \\
     \langle \psi_{\ell_a}, \mathcal{B} e_1\rangle_H & \ldots  & \langle \psi_{\ell_a}, \mathcal{B} e_m\rangle_H \end{bmatrix} \in \mathbb{R}^{\ell_a \times m}\;.
\end{equation*}
Thus, collectively, the POD-G approximation of $\mathrm{OC}_{\scaleto{\mathrm{FOM}}{4pt}}$ is given as
\begin{subnumcases}{\mathrm{OC}_{\scaleto{\mathrm{POD-G}}{4pt}}:=}
     \dot{\alpha}(t) = A_{\ell} \alpha(t) + B_{\ell} u(t)  \label{eq_def:OC_FOTR_PODG_1},  \\
    \alpha(0)  = \alpha_0,\label{eq_def:OC_FOTR_PODG_2}\\[1em]
    -\dot{\lambda}^{\ell_a}(t) = A_{\ell_a} \lambda^{\ell_a}(t) + M \alpha(t) - \hat{q}_\mathrm{d}(t)\label{eq_def:OC_FOTR_PODG_3},  \\
    \lambda^{\ell_a}(T) = 0, \label{eq_def:OC_FOTR_PODG_4}\\[1em]
    \mu u(t) + B_{\ell_a}^\top \lambda^{\ell_a}(t) = 0 \, .\label{eq_def:OC_FOTR_PODG_5}
\end{subnumcases}
We also note that the cost functional for the POD-G approximation of $\mathrm{OC}_{\scaleto{\mathrm{FOM}}{4pt}}$  
\begin{equation}\label{eq:PODG_costFunc_continuous}
    \mathcal{J}_{\scaleto{\mathrm{POD-G}}{4pt}}(\alpha, u) = \frac{1}{2} \int^{T}_0 \norm{\sum^{\ell}_{i=1} \alpha_i(t)\phi_i - q_\mathrm{d}(t)}^2_H \, \mathrm{d}t + \frac{\mu}{2} \int_0^{T}\norm{u(t)}^2_{U}\,\mathrm{d}t
\end{equation}
is now used as an approximation of \eqref{eq:FOM_costFunc_continuous}.

\subsection{sPOD-Galerkin approximation of $\mathrm{OC}_{\scaleto{\mathrm{FOM}}{4pt}}$}
Given the trajectory $q \in C([0, T]; \mathcal{D}(\mathcal{A}))$ from the FOM state equation, the sPOD decomposes it with a nonlinear ansatz,
\begin{equation}\label{eq:sPOD_advection}
    q(t) \approx \sum^{\tilde{\ell}}_{i=1} \alpha_{i}(t) \mathcal{T}(z(t)) \phi_{i}\;,
\end{equation}
where $z(t) \in \mathbb{R}$ is the single time-dependent shift, $\mathcal{T}\colon \mathbb{R} \rightarrow \mathcal{L}(H)$ is the corresponding transformation operator, $\alpha_i \in L^2(0, T)$ are the reduced state variables and $\{\phi_i\}^{\tilde{\ell}}_{i=1} \in \mathcal{D}(\mathcal{A})$ are stationary modes.
\begin{rmrk}
    The choice of the sPOD ansatz \eqref{eq:sPOD_advection} used is motivated by the dynamics of the linear advection equation considered here with a single traveling wave. 
A more general version of the ansatz is discussed in \cite[Eq.~1.4]{black_projection-based_2020}.
\end{rmrk}
To minimize the approximation error in \eqref{eq:sPOD_advection}, we solve the following minimization problem
\begin{subnumcases}{} \label{eq:sPOD_advection_optimization}
    \mathrm{min} \: \frac{1}{2} \int^{T}_0 \norm{q(t) - \mathcal{T}(z(t))\sum^{\tilde{\ell}}_{i=1} \alpha_i(t) \phi_i }^2_H \mathrm{d}t ,\\[1em]
    \text{s.t.} \: \{\phi_i\}^{\tilde{\ell}}_{i=1} \subset \mathcal{D}(\mathcal{A}) \: \text{and} \: \langle \phi_i, \phi_j \rangle_H = \delta_{ij}, \alpha_i \in L^2(0, T) \;
\end{subnumcases}
with the assumption that the shift function $z\in L^2(0,T)$ is given. 
This minimization problem can also be viewed as a POD minimization of the shifted data $\mathcal{T}^*(z(t))q(t), t\in [0,T]$. 
Also, keeping in mind the dynamics of the problem at hand, the transformation operators are taken to be periodic shift operators and are given as 
\begin{align*}
    \mathcal{T}(z)\phi(x)=
    \begin{cases}
        \phi(x-\eta) & \text{for } \eta \le x \le l, \\
        \phi(x-\eta+l) & \text{for } 0 \le x < \eta 
    \end{cases}
\end{align*}
with $\eta\vcentcolon=z\mod l$, see e.g.~\cite[Def.~1.2.2]{schulze_energy-based}.
Such a choice results in $\mathcal{T}$ being isometric and $H^1_\mathrm{per}(0, l)$ being a $\mathcal{T}(z)$-invariant subspace for any $z\in \mathbb{R}$.
\begin{rmrk}\label{rem:diff_T_vs_phi}
 If $\phi$ is sufficiently regular, we have $\mathcal{T}'(z)\phi:=\tfrac{\mathrm{d}}{\mathrm{d}z}(\mathcal{T}(z)\phi)=-\mathcal{T}(z)\phi'$ and $\mathcal{T}''(z)\phi=\mathcal{T}(z)\phi''$. 
Taking into account the isometry of $\mathcal{T}$, it holds that $\mathcal{T}^*(z)=\mathcal{T}(-z)=\mathcal{T}^{-1}(z).$
\end{rmrk}
Subsequently, the sPOD-G ROM for the FOM state equation is obtained by a Galerkin projection and is given as
\begin{equation}\label{eq:sPOD_galerkin_advection}
    \begin{bmatrix} I_{\tilde{\ell}} & N \alpha(t) \\ \alpha(t)^\top N^\top & \alpha(t)^\top M_2 \alpha(t)\end{bmatrix} \begin{bmatrix} \dot{\alpha}(t) \\ \dot{z}(t) \end{bmatrix} = v \begin{bmatrix} N & 0 \\ \alpha(t)^\top M_2 & 0 \end{bmatrix}
    \begin{bmatrix} \alpha(t) \\ z(t) \end{bmatrix} + \begin{bmatrix}  B_1(z(t)) \\ \alpha(t)^\top B_2(z(t))\end{bmatrix} u(t)
\end{equation}
which is a nonlinear equation where
\begin{align*}
N &:= -\big[\langle \phi_i,\;\phi_j'\rangle_H\big]_{i,j=1}^{\tilde{\ell}},
&M_2 &:= \big[\langle \phi_i',\;\phi_j'\rangle_H\big]_{i,j=1}^{\tilde{\ell}},\\
    B_1(z(t)) &:= \big[\langle \mathcal{T}(z(t))\phi_i,\;\mathcal{B}e_j\rangle_H\big]_{i,j=1}^{\tilde{\ell},m}, &B_2(z(t)) &:= \big[\langle \mathcal{T}'(z(t))\phi_i,\;\mathcal{B}e_j\rangle_H\big]_{i,j=1}^{\tilde{\ell},m}.
\end{align*}
For given $z(0)$, the initial value $\alpha(0)$ for the sPOD-G ROM state equation is chosen so that the approximation error in 
\begin{equation*}
    q_0 \approx \sum^{\tilde{\ell}}_{i=1} \alpha_i(0) \mathcal{T}(z(0)) \phi_i
\end{equation*}
is minimized. 
Performing a Galerkin projection of the above expression onto the shifted modes $\{\mathcal{T}(z(0)) \phi_j\}^{\tilde{\ell}}_{j=1}$ leads to
\begin{align*}
    \langle \mathcal{T}(z(0)) \phi_j, q_0\rangle_H &\approx \Bigl \langle \mathcal{T}(z(0)) \phi_j, \sum^{\tilde{\ell}}_{i=1} \alpha_i(0) \mathcal{T}(z(0)) \phi_i \Bigr\rangle_H  \nonumber
    = \Bigl \langle \phi_j, \sum^{\tilde{\ell}}_{i=1} \alpha_i(0) \phi_i \Bigr\rangle_H \nonumber.
\end{align*}
We additionally assume $z(0) = 0$, which results in $\mathcal{T}(z(0)) = I_n$,
therefore, it follows that $\alpha(0) \approx [\langle \phi_j, q_0\rangle_H]^{\tilde{\ell}}_{j=1}$ corresponds to the orthogonal projection of the initial value $q_0$ onto the span of the dominant stationary modes $\{\phi_j\}^{\tilde{\ell}}_{j=1}$.

For the reduced-order approximation of the adjoint equation of the FOM \eqref{eq:OC_FOM_3}, \eqref{eq:OC_FOM_4}, we follow similar steps as for the state equation and immediately formulate the sPOD ansatz for the adjoint $\lambda \in C([0, T];\mathcal{D}(\mathcal{A}^*))$ as
\begin{equation}\label{eq:sPODG_FOTR_adjoint_ansatz}
    \lambda(t) = \sum^{\tilde{\ell}_a}_{i=1} \lambda^{\ell_a}_i(t) \mathcal{T}(z^{\ell_a}(t))\psi_i, \; z^{\ell_a}(t) \in \mathbb{R}, \; \lambda^{\ell_a}_i \in L^2(0, T), \; \psi_i \in D(\mathcal{A}^*),\; \langle \psi_i,\psi_j\rangle_H=\delta_{ij}\; \text{for }i,j=1,\ldots,\tilde{\ell}_a.
\end{equation}
The adjoint is generally approximated by using a different shift $z^{\ell_a}$ than for the state. This is due to the difference in the propagation velocity of the adjoint from that of the FOM state in the presence of an arbitrary control term. 
Now, since the FOM adjoint equation is linear and $\mathcal{A}^*$  is also the generator of a strongly continuous (semi-)group, we follow similar arguments as above and arrive at the following sPOD-G ROM for the adjoint
\begin{equation}
        \begin{bmatrix} I_{\tilde{\ell}_a} & \tilde{N} \lambda^{\ell_a} \\ (\lambda^{\ell_a})^\top \tilde{N}^\top & (\lambda^{\ell_a})^\top \tilde{M}_2 \lambda^{\ell_a}\end{bmatrix} \begin{bmatrix} \dot{\lambda}^{\ell_a} \\ \dot{z}^{\ell_a} \end{bmatrix} = - v\begin{bmatrix} \tilde{N} & 0 \\ (\lambda^{\ell_a})^\top \tilde{M}_2 & 0 \end{bmatrix}
    \begin{bmatrix} \lambda^{\ell_a} \\ z^{\ell_a} \end{bmatrix} 
    - \begin{bmatrix}  \tilde{M}(z^{\ell_a}- z)\alpha - \tilde{q}_\mathrm{d}(t,z^{\ell_a}) \\ (\lambda^{\ell_a})^\top(\doubletilde{M}(z^{\ell_a}- z)\alpha - \doubletilde{q}_\mathrm{d}(t,z^{\ell_a}))\end{bmatrix}
\end{equation}
with 
\begin{align*}
\tilde{N} &:= -\big[\langle \psi_i,\;\psi_j'\rangle_H\big]_{i,j=1}^{\tilde{\ell}_a},
&\tilde{M}_2 &:= \big[\langle \psi_i',\;\psi_j'\rangle_H\big]_{i,j=1}^{\tilde{\ell}_a},
\end{align*}
\begin{equation*}
    \tilde{M}(\eta) := \begin{bmatrix} \langle \mathcal{T}(\eta)\psi_1, \phi_1\rangle_H  & \ldots  & \langle \mathcal{T}(\eta)\psi_1, \phi_{\tilde{\ell}}\rangle_H \\ \vdots & \ddots & \vdots \\
     \langle \mathcal{T}(\eta)\psi_{\tilde{\ell}_a}, \phi_1\rangle_H & \ldots  & \langle \mathcal{T}(\eta)\psi_{\tilde{\ell}_a}, \phi_{\tilde{\ell}}\rangle_H \end{bmatrix} \in \mathbb{R}^{\tilde{\ell}_a \times \tilde{\ell}}\;,
\end{equation*}
\begin{equation*}
    \doubletilde{M}(\eta) := \begin{bmatrix} \langle \mathcal{T}'(\eta)\psi_1, \phi_1\rangle_H  & \ldots  & \langle \mathcal{T}'(\eta)\psi_1, \phi_{\tilde{\ell}}\rangle_H \\ \vdots & \ddots & \vdots \\
     \langle \mathcal{T}'(\eta)\psi_{\tilde{\ell}_a}, \phi_1\rangle_H & \ldots  & \langle \mathcal{T}'(\eta)\psi_{\tilde{\ell}_a}, \phi_{\tilde{\ell}}\rangle_H \end{bmatrix} \in \mathbb{R}^{\tilde{\ell}_a \times \tilde{\ell}}\;,
\end{equation*}
and $\tilde{q}_\mathrm{d}(t,z^{\ell_a}) = [\langle \mathcal{T}(z^{\ell_a})\psi_j, q_\mathrm{d}(t)\rangle_H ]^{\tilde{\ell}_a}_{j=1} \in \mathbb{R}^{\tilde{\ell}_a}$, $\doubletilde{q}_\mathrm{d}(t,z^{\ell_a}) = [\langle \mathcal{T}'(z^{\ell_a})\psi_j, q_\mathrm{d}(t)\rangle_H ]^{\tilde{\ell}_a}_{j=1} \in \mathbb{R}^{\tilde{\ell}_a}$.
We note that for the derivation of the terms involving $\tilde{M}$ and $\doubletilde{M}$, we exploited Remark~\ref{rem:diff_T_vs_phi} and the fact that $\mathcal{T}$ is a group homomorphism, see e.g.~\cite[App.~A]{schulze_energy-based}, to obtain $\mathcal{T}(z)^*\mathcal{T}(z^{\ell_a}) = \mathcal{T}(z^{\ell_a}-z)$.
Since the sPOD-G adjoint is solved backward in time, the terminal conditions $z^{\ell_a}(T)$ and $\lambda^{\ell_a}(T)$ need to be prescribed and be compatible with the terminal condition for the full-order adjoint state, i.e.,
\begin{equation*}
    \sum^{\tilde{\ell}_a}_{i=1} \lambda_i^{\ell_a}(T) \mathcal{T}(z^{\ell_a}(T)) \psi_i = \lambda(T) = 0\;.
\end{equation*}
Since $\mathcal{T}$ is pointwise isometric and the modes $\psi_1,\ldots,\psi_{\tilde{\ell}_a}$ are orthonormal, this equation is equivalent to $\lambda_1^{\ell_a}(T)=\ldots=\lambda_{\tilde{\ell}_a}^{\ell_a}(T)=0$.
We note that this terminal condition is satisfied for any shift value $z^{\ell_a}(T)$, see also \cite[Remark 5.8]{black_projection-based_2020}.

Lastly, the optimality condition on the control can be approximated as
\begin{equation}
    \mu u(t) + \tilde{B}(z^{\ell_a}(t))^\top\lambda^{\ell_a}(t) = 0
\end{equation}
where
\begin{equation*}
    \tilde{B}(z^{\ell_a}(t)) := \begin{bmatrix} \langle \mathcal{T}(z^{\ell_a}(t))\psi_1, \mathcal{B}e_1\rangle_H  & \ldots  & \langle \mathcal{T}(z^{\ell_a}(t))\psi_{1}, \mathcal{B}e_m\rangle_H \\ \vdots & \ddots & \vdots\\
     \langle \mathcal{T}(z^{\ell_a}(t))\psi_{\tilde{\ell}_a}, \mathcal{B} e_1 \rangle_H & \ldots  & \langle \mathcal{T}(z^{\ell_a}(t))\psi_{\tilde{\ell}_a}, \mathcal{B} e_m\rangle_H \end{bmatrix} \in \mathbb{R}^{\tilde{\ell}_a \times m }\;.
\end{equation*}
Thus, collectively, the sPOD-G approximation of $\mathrm{OC}_{\scaleto{\mathrm{FOM}}{4pt}}$ is given as
\begin{subnumcases}{\mathrm{OC}_{\scaleto{\mathrm{sPOD-G}}{4pt}}:=}
    \begin{bmatrix} I_{\tilde{\ell}} & N \alpha(t) \\ \alpha(t)^\top N^\top & \alpha(t)^\top M_2 \alpha(t)\end{bmatrix} \begin{bmatrix} \dot{\alpha}(t) \\ \dot{z}(t) \end{bmatrix} = v \begin{bmatrix} N & 0 \\ \alpha(t)^\top M_2 & 0 \end{bmatrix}
    \begin{bmatrix} \alpha(t) \\ z(t) \end{bmatrix} + \begin{bmatrix}  B_1(z(t)) \\ \alpha(t)^\top B_2(z(t))\end{bmatrix} u(t)\label{eq:OC_FOTR_sPODG_1},  \\
    \alpha(0) = \alpha_0, \quad \quad z(0) = 0 \label{eq:OC_FOTR_sPODG_2},\\[1em]
    \begin{split}
        \begin{bmatrix} I_{\tilde{\ell}_a} & \tilde{N} \lambda^{\ell_a}(t) \\ \lambda^{\ell_a}(t)^\top \tilde{N}^\top & \lambda^{\ell_a}(t)^\top \tilde{M}_2 \lambda^{\ell_a}(t)\end{bmatrix} \begin{bmatrix} \dot{\lambda}^{\ell_a}(t) \\ \dot{z}^{\ell_a}(t) \end{bmatrix} = - v\begin{bmatrix} \tilde{N} & 0 \\ \lambda^{\ell_a}(t)^\top \tilde{M}_2 & 0 \end{bmatrix}
    \begin{bmatrix} \lambda^{\ell_a}(t) \\ z^{\ell_a}(t) \end{bmatrix} 
    \\ - \begin{bmatrix}  \tilde{M}(z^{\ell_a}(t)- z(t))\alpha(t) - \tilde{q}_\mathrm{d}(t,z^{\ell_a}(t)) \\ \lambda^{\ell_a}(t)^\top(\doubletilde{M}(z^{\ell_a}(t)-z(t))\alpha(t) - \doubletilde{q}_\mathrm{d}(t,z^{\ell_a}(t)))\end{bmatrix}
    \end{split} \label{eq:OC_FOTR_sPODG_3},  \\
    \lambda^{\ell_a}(T) = 0  \text{ and} \; z^{\ell_a}(T) \in \mathbb{R} \label{eq:OC_FOTR_sPODG_4}\\[1em] 
    \mu u(t) + \tilde{B}(z^{\ell_a}(t))^\top\lambda^{\ell_a}(t) = 0. \label{eq:OC_FOTR_sPODG_5}
\end{subnumcases}
We note that the sPOD-G approximation is parameterized by the terminal shift value $z^{\ell_a}(T)$.
Analogously to the POD-G case, the cost functional $\mathcal{J}$ from \eqref{eq:FOM_costFunc_continuous} is approximated as
\begin{equation}\label{eq_def:sPOD_costFunc_continuous}
\mathcal{J}_{\scaleto{\mathrm{sPOD-G}}{4pt}}(\alpha, z, u):= \frac{1}{2} \int^{T}_0 \norm{\sum^{\tilde{\ell}}_{i=1} \alpha_i(t)\mathcal{T}(z(t))\phi_i - q_\mathrm{d}(t)}^2_H \, \mathrm{d}t + \frac{\mu}{2} \int_0^{T}\norm{u(t)}^2_{U}\,\mathrm{d}t\;.
\end{equation}

\begin{rmrk} 
Reduced optimality conditions can also be formulated based on the first-reduce-then-optimize (FRTO) framework, see \cite[Sec. 3]{breiten2026optimalcontrolshiftedproper}. 
    For POD-G, considering linear systems with $\phi=\psi$ the FOTR and FRTO frameworks leads to the same optimality system.
    A detailed treatment of both frameworks in the POD-G context can be seen in \cite{negri_reduced_2013, negri_reduced_2015}.
    However, for the sPOD-G method, $\mathrm{OC}_{\scaleto{\mathrm{sPOD-G}}{4pt}}$ and the optimality conditions of FRTO from \cite[Eq. 36(a)-(e)]{breiten2026optimalcontrolshiftedproper} do not appear to commute. 
\end{rmrk}

\section{Algorithmic details}
In this section, we explain the algorithmic aspects of using ROMs in the optimal control problem in an FOTR framework. 
Since POD-G based optimal control problems and corresponding algorithms to solve such problems have been widely studied, we refrain ourselves from specifically outlining the algorithm here. For some of such algorithms, see \cite{Afanasiev_wake_2001, ravindran_adaptive_2002}. 
Firstly, the FOM state and the adjoint equations are solved using a first-order upwind scheme. 
Next, for the ROMs, we use an explicit Euler scheme for both the POD-G and sPOD-G reduced-state equations. 
For the ROM adjoint equations for both methods, we again use an explicit Euler scheme.

\begin{rmrk}
POD-G results in a linear reduced-state and reduced-adjoint equation.
However, sPOD-G results in a nonlinear reduced-state as well as a nonlinear reduced-adjoint equation.
    For our example problem, using an explicit first-order time integrator suffices for both these nonlinear equations.
    However, this might not always be the case and the sPOD-G based FOTR framework could end up needing higher-order time integrators for stability in more complex scenarios.
\end{rmrk}
Algorithm \ref{alg:FOTR_sPODG} summarizes the required steps for using the sPOD-G method in our optimal control context considered here.
\begin{algorithm}  
  \caption{Optimal control with sPOD-G}\label{alg:FOTR_sPODG}
  \begin{algorithmic}[1]
  \Require{$q_0$, $q_\mathrm{d}$, $\tilde{\ell}$, $\tilde{\ell}_a$, $\mu$, $\omega^0$, $n_\mathrm{iter}$, $\beta$, $n_\mathrm{samples}$, $n_\mathrm{samples}^{\ell_a}$} 
    \State \textbf{Initialize:} $u=u_0$
      \For{$i = 1, \ldots, n_{\mathrm{iter}}$}
        \If{$\mathrm{refine}$} 
            \State $q= \textsc{State}(u^i, q_0)$ \Comment{Solve \eqref{eq:OC_FOM_1} and \eqref{eq:OC_FOM_2}}
            \State $\lambda= \textsc{Adjoint}(q, q_\mathrm{d})$ \Comment{Solve \eqref{eq:OC_FOM_3} and \eqref{eq:OC_FOM_4}}
            \State $\{\mathcal{T}(z)\phi_k\}^{\tilde{\ell}}_{k=1} = \textsc{Basis} ([q(t_1) \: \ldots \: q(t_{n_\mathrm{t}})], \tilde{\ell}, n_\mathrm{samples})$
            \State $\{\mathcal{T}(z^{\ell_a})\psi_k\}^{\tilde{\ell}_a}_{k=1} = \textsc{Basis} ([\lambda(t_1) \: \ldots \: \lambda(t_{n_\mathrm{t}})], \tilde{\ell}_a, n_\mathrm{samples}^{\ell_a})$
        \EndIf
        \State $\alpha^i, z^i = \textsc{ReducedState}(\{\mathcal{T}(z)\phi_k\}^{\tilde{\ell}}_{k=1}, u^i, q_0)$ \Comment{Solve \eqref{eq:OC_FOTR_sPODG_1} and \eqref{eq:OC_FOTR_sPODG_2}}
        \State $\lambda^{\ell_a, i}, z^{\ell_a, i} = \textsc{ReducedAdjoint}(\{\mathcal{T}(z^{\ell_a})\psi_k\}^{\tilde{\ell}_a}_{k=1}, q_\mathrm{d})$ \Comment{Solve \eqref{eq:OC_FOTR_sPODG_3} and \eqref{eq:OC_FOTR_sPODG_4}}
        \State $\tfrac{\mathrm{d}\mathcal{L}}{\mathrm{d}u^i}=\textsc{Gradient}(\lambda^{\ell_a, i}, z^{\ell_a, i}, \{\mathcal{T}(z^{\ell_a})\psi_k\}^{\tilde{\ell}_a}_{k=1}, u^i)$ \Comment{Solve \eqref{eq:OC_FOTR_sPODG_5}}
        \State $\omega^i = \textsc{StepSize}(\omega^{i- 1}, \tfrac{\mathrm{d}\mathcal{L}}{\mathrm{d} u^i}, \{\mathcal{T}(z)\phi_k\}^{\tilde{\ell}}_{k=1}, \alpha^i, u^{i})$
        \State $u^{i+1} = u^{i} - \omega^i\left(\tfrac{\mathrm{d}\mathcal{L}}{\mathrm{d} u^i}\right)$ 
        \If{$i==n_{\mathrm{iter}}$} 
            \State break
        \ElsIf{$\norm{\tfrac{\mathrm{d}\mathcal{L}}{\mathrm{d} u^i}} / \norm{\tfrac{\mathrm{d}\mathcal{L}}{\mathrm{d} u^1}} < \beta$}
            \State set $u=u^{i+1}$ and return
        \EndIf
      \EndFor
  \Ensure{$u$}
  \end{algorithmic}
\end{algorithm}
Since the reduced-state equations \eqref{eq:OC_FOTR_sPODG_1}, \eqref{eq:OC_FOTR_sPODG_2} and the reduced-adjoint equations \eqref{eq:OC_FOTR_sPODG_3}, \eqref{eq:OC_FOTR_sPODG_4} in the above algorithm are nonlinear, the corresponding nonlinear quantities must be assembled at every time step.
However, repeatedly forming these terms is computationally expensive, particularly when they do not scale with the reduced-order dimension.
To address this issue, we precompute most terms that are nonlinear and depend on the FOM dimension.
In particular, we evaluate the state-shift-dependent terms $B_1(z(t))$, $B_2(z(t))$ and the adjoint-shift-dependent term $\tilde{B}(z^{{\ell}_a}(t))$ at $n_{\mathrm{samples}}$ and $n_{\mathrm{samples}}^{\ell_a}$ shift samples for the state and adjoint, respectively.
For terms containing both the state and adjoint shifts, $\tilde{M}(z^{\ell_a}(t) - z(t))$, $\doubletilde{M}(z^{{\ell}_a}(t)- z(t))$, we effectively do not sample anything for our implementation since $z^{\ell_a} - z = 0$ (see Rem.~\ref{rem:adjoint_shift}).
During time integration, we then use linear interpolation to recover the values of $B_1(z(t))$, $B_2(z(t))$ and $\tilde{B}(z^{{\ell}_a}(t))$.
Consequently, the expensive construction of these terms is shifted to a preprocessing stage prior to the actual time integration.
This strategy is analogous to the approach used in \cite{black_efficient_2021}.
We note that we cannot use the same strategy for $\tilde{q}_\mathrm{d}$ and $\doubletilde{q}_\mathrm{d}$ since they depend in addition on the time-dependent target state $q_\mathrm{d}$.
Thus, the evaluation of these two terms in the online phase still scales with the dimension of the (semi-discretized) FOM.

During the refinement steps shown in algorithm \ref{alg:FOTR_sPODG}, the sPOD ansatz is applied to the state and adjoint FOM snapshots to construct the shifted basis.
 During construction of the shifted basis, the shifts $z$ and $z^{\ell_a}$ and their corresponding transformation operators also need to be constructed.
 Although these operators are sparse, we need to construct $n_t$ of these operators for both the state and the adjoint, making it computationally costly.
 To reduce this expense, we compute the shifts via a template fitting approach \cite{rowley_reconstruction_2000} and their corresponding transformation operators once from the uncontrolled profile and keep them fixed throughout the optimization. 
However, $\{\phi_k\}$ and $\{\psi_k\}$ need to be updated and are obtained by performing an SVD on the shifted state  and the shifted adjoint snapshot matrices, respectively.
The basis refinement is performed at every fifth optimization step, as well as whenever the step-size selection criterion fails.
For the step size selection, we use a combination of two-way backtracking \cite{truong_backtracking_2021} and the Barzilai-Borwein (BB) method \cite{barzilai_two-point_1988}. 
During the initial steps where $\norm{\tfrac{\mathrm{d}\mathcal{L}}{\mathrm{d} u^i}} / \norm{\tfrac{\mathrm{d}\mathcal{L}}{\mathrm{d} u^1}} > 5 \times 10^{-3}$, the two-way backtracking is used to select the step size, whereas when the value falls below $5 \times 10^{-3}$, the BB method is used. 
\begin{rmrk}\label{rem:adjoint_shift}
    Within the sPOD-G ansatz of the adjoint equation, as outlined in \eqref{eq:sPODG_FOTR_adjoint_ansatz}, it becomes evident that the adjoint equation requires its unique shift value. 
    However, under certain idealized conditions, such as those considered in our experimental scenario, it is possible to use the same shift value to transform both the state and adjoint into the stationary frame. 
    This assumption is built on the notion that the transport dynamics is predominantly influenced by the advection operator and that the control input does not significantly alter the overall propagation velocity of the advected pulse. 
    Given this presumption, it is sufficient to perform the computation of transformation operators solely for the state equation, which can then be subsequently applied to the adjoint system as well. 
    Alternatively, one could choose to independently calculate the shift for the adjoint system, but this would require an additional computation of the shift and the assembly of transformation operators, resulting in an increased computational overhead.
    Moreover, in that case it would be not clear how to actually choose the terminal shift for the adjoint equation.
\end{rmrk}
\begin{rmrk}
    We observe that for \eqref{eq:OC_FOTR_sPODG_3}, at $t=T$ the coefficient matrix of the left-hand side is given by
    \begin{equation*}
        \begin{bmatrix} I_{\tilde{\ell}_a} & \tilde{N} \lambda^{\ell_a}(T) \\ \lambda^{\ell_a}(T)^\top \tilde{N}^\top & \lambda^{\ell_a}(T)^\top \tilde{M}_2 \lambda^{\ell_a}(T)\end{bmatrix} = \begin{bmatrix} I_{\tilde{\ell}_a} & 0 \\ 0 & 0\end{bmatrix},
    \end{equation*}
    which is singular.
    To make the problem numerically solvable, we regularize the coefficient matrix by adding a small $\epsilon = 10^{-14}$ to the last entry of the diagonal.
\end{rmrk}

\section{Numerical results}
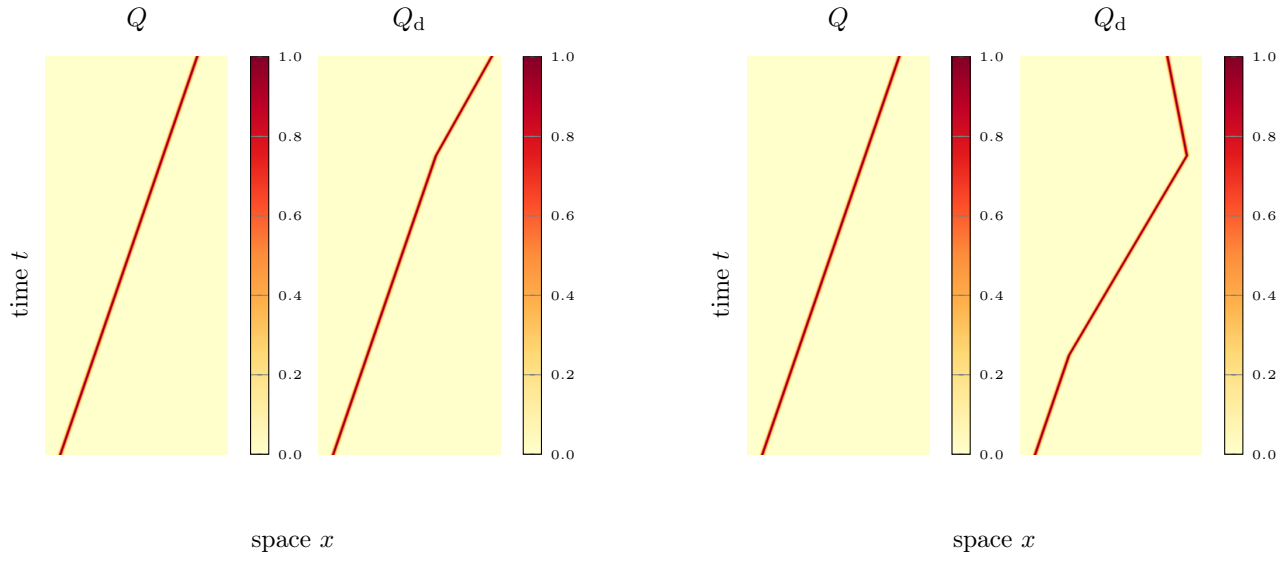
\begin{figure}[htp!]
    \centering
    \begin{subfigure}[t]{0.45\textwidth}
        \centering
        \setlength\figureheight{0.7\linewidth}
        \setlength\figurewidth{0.7\linewidth}
        \input{State_target_p1} % Input the first plot
        \caption{Single tilt: Traveling wave with kink at $t=\frac{3}{4}T$.}
        \label{fig:Example_1}
    \end{subfigure}
    \hspace{0.05\textwidth}
    \begin{subfigure}[t]{0.45\textwidth}
        \centering
        \setlength\figureheight{0.7\linewidth}
        \setlength\figurewidth{0.7\linewidth}
        \input{State_target_p4} % Input the second plot
        \caption{Double tilt: Traveling wave with kinks at $t=\frac{1}{4}T$ and $t=\frac{3}{4}T$.}
        \label{fig:Example_4}
    \end{subfigure}
    \caption{Plots for the state and the target for the two example problems}
    \label{fig:advection_and_target}
\end{figure}

\begin{figure}[htp!]
    \centering
    \includegraphics[scale=0.42]{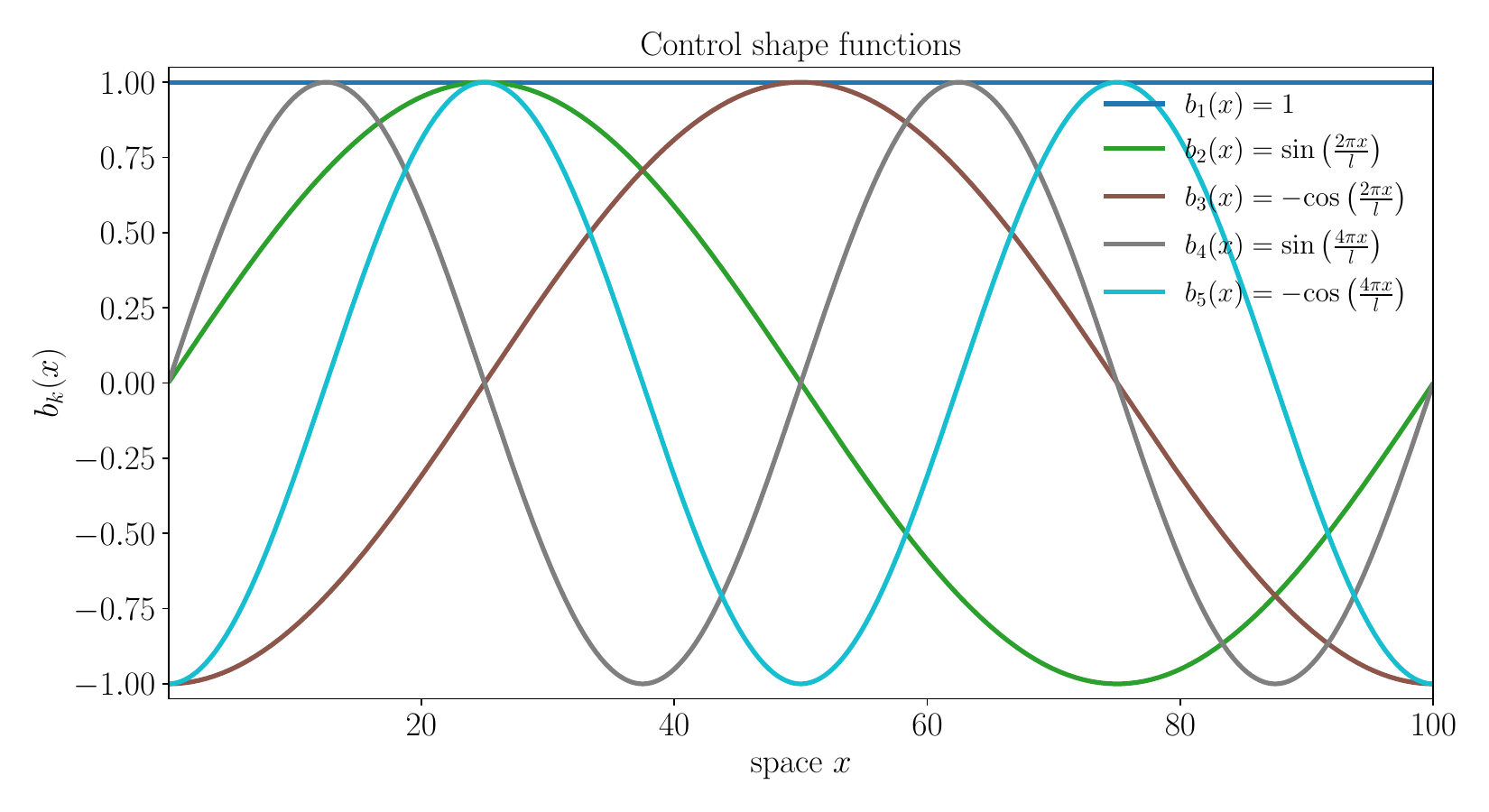}
    \caption{Control shape functions}
    \label{fig:Controls}
\end{figure}

\begin{figure}[htp!]
    \centering
        \begin{subfigure}[t]{0.48\textwidth}
        \centering
        \setlength\figureheight{0.8\linewidth}
        \setlength\figurewidth{0.9\linewidth}
        \input{FOTR_modesVsCost_Adaptive_P1}   % Input the second plot
        \caption{Single tilt problem}
        \label{fig:J_vs_modes_FOTR_AB_shifting}
    \end{subfigure}
    \hspace{0.01\textwidth}
    \begin{subfigure}[t]{0.48\textwidth}
        \centering
        \setlength\figureheight{0.8\linewidth}
        \setlength\figurewidth{0.9\linewidth}
        \input{FOTR_modesVsCost_Adaptive_P}   % Input the second plot
        \caption{Double tilt problem}
        \label{fig:J_vs_modes_FOTR_AB_shifting_3}
    \end{subfigure}
    \caption{Plots for $\mathcal{J}$ vs.~$\mathrm{modes}$ for both example problems}
    \label{fig:J_shifting}
\end{figure}
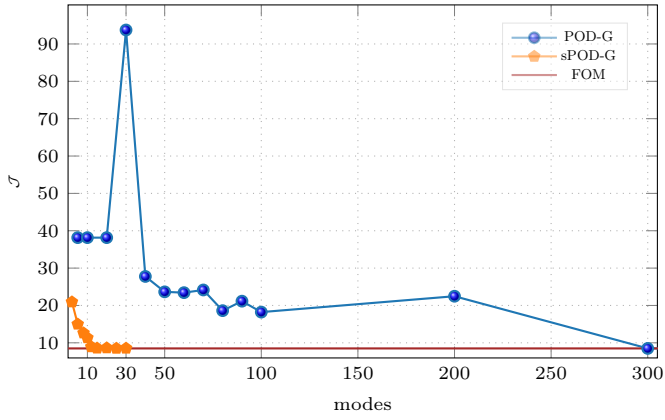
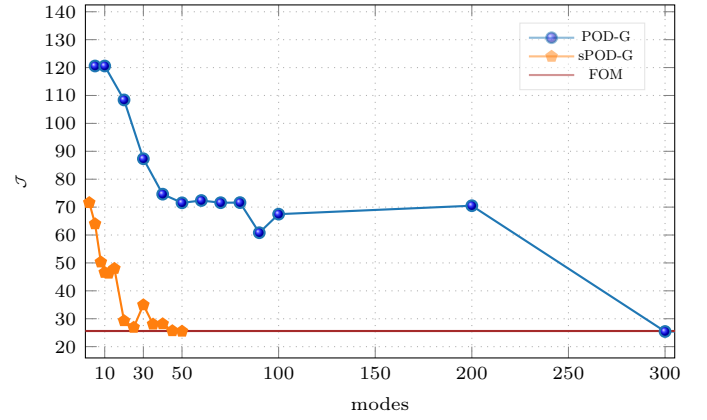

All numerical tests \footnote{source code can be found in \urlstyle{tt}\url{https://doi.org/10.5281/zenodo.20402350}} were run using Python 3.12 on a Macbook Air M1(2020) with an 8-core CPU and 16 GB of RAM. 
We consider two different test cases as shown in Figure~\ref{fig:advection_and_target} where $Q$ and $Q_\mathrm{d}$ are snapshot matrices constructed by stacking the FOM state and the target state column-wise, respectively.
To obtain the solution, a one-dimensional strip of length $l=100$ is discretized into grid points $n=3201$ in the spatial domain with $\Delta x = 0.03124$ with periodic boundary conditions.
The PDE is solved with the following initial condition
\begin{equation}\label{eq_def:advection_init}
    q_0(x) := \mathrm{exp} \; \left(-\left(x - \frac{l}{12}\right)^2\right) \;
\end{equation}
until the final time $T=136.2642$.
The time domain is discretized into $n_t = 2400$ time steps with $\Delta t= 0.0568$.
The advection velocity is $v=0.55$ and the starting control $u = 0$.
The control shape functions in our tests are chosen as
\begin{equation}\label{eq_def:mask_sin_cos}
    b_1(x) = 1, \quad b_{2k}(x) = \mathrm{sin}\left(\frac{2\pi k x}{l}\right), \quad b_{2k + 1}(x) = -\mathrm{cos}\left(\frac{2\pi k x}{l}\right)
\end{equation}
for $k=1, \ldots, \xi$, where $\xi = 20$  and $m=2 \xi + 1 = 41$. 
The first five are shown in Figure \ref{fig:Controls}.
\begin{rmrk}\label{rem:control_shape}
    The choice of $b_k$ plays a decisive role in the overall convergence of the optimal control problem for sPOD-G.
    In our tests, $b_k$ from \eqref{eq_def:mask_sin_cos} are essentially the real vectors representing pairs of conjugate complex eigenfunctions of the operator $\mathcal{A} = -v \frac{\mathrm{d}}{\mathrm{d}x}$.
    Considering such control shape functions provides an upper bound on the rank of the snapshot matrix in the stationary frame, see \cite[Prop 2.2]{breiten2026optimalcontrolshiftedproper}
\end{rmrk}
The optimization parameters used for the example problems are given in Table \ref{tab:params_1}.
\begin{table}[h!]
\begin{center}
\begin{minipage}{\textwidth}
\centering
\begin{tabular}{|l||c|c|c|c|c|c|c|} 
 \hline
 & $\mu$ & $\beta$ & $\omega^0$ & $n_{\mathrm{iter}}$ & $n_{\mathrm{samples}}$ & $n^{\ell_a}_{\mathrm{samples}}$\\
\hline 
Single tilt problem & $10^{-3}$ & $1\times 10^{-5}$ & $1$ & $20000$ & $800$ & $800$\\
\hline   
Double tilt problem & $10^{-3}$ & $1\times 10^{-5}$ & $1$ & $20000$ & $800$ & $800$\\
\hline  
\end{tabular}
\caption{Constant parameters}
\label{tab:params_1}
\end{minipage}
\end{center}
\end{table}
Subsequently, the results for both example problems are shown in Figure~\ref{fig:J_shifting}.
For the single tilt problem, cf.~Figure~\ref{fig:J_vs_modes_FOTR_AB_shifting}, the sPOD-G method needs close to $25$ modes to reach the FOM cost functional, whereas POD-G requires close to $300$ modes for the same.
Regarding the double tilt problem, Figure~\ref{fig:J_vs_modes_FOTR_AB_shifting_3} shows that sPOD-G needs close to $45$ modes to reach the FOM cost, whereas POD-G needs $300$ modes for the same.
A point to note here is that for the mode-based study shown in Figure~\ref{fig:J_shifting}, the number of truncated modes for the state and the adjoint equations is chosen to be equal for simplicity.
\begin{figure}[htp!]
    \centering
        \begin{subfigure}[t]{0.48\textwidth}
        \centering
        \setlength\figureheight{0.8\linewidth}
        \setlength\figurewidth{0.9\linewidth}
        \input{FOTR_tolVsCost_Adaptive_P1}   % Input the second plot
        \caption{Single tilt problem}
        \label{fig:J_vs_tol_FOTR_AB_shifting}
    \end{subfigure}
    \hspace{0.01\textwidth}
    \begin{subfigure}[t]{0.48\textwidth}
        \centering
        \setlength\figureheight{0.8\linewidth}
        \setlength\figurewidth{0.9\linewidth}
        \input{FOTR_tolVsCost_Adaptive_P}   % Input the second plot
        \caption{Double tilt problem}
        \label{fig:J_vs_tol_FOTR_AB_shifting_3}
    \end{subfigure}
    \caption{Plots for $\mathcal{J}$ vs.~ $\mathrm{tol}$ for both example problems with the respective state and adjoint truncated modes for each tolerance value.}
    \label{fig:J_tol_shifting}
\end{figure}
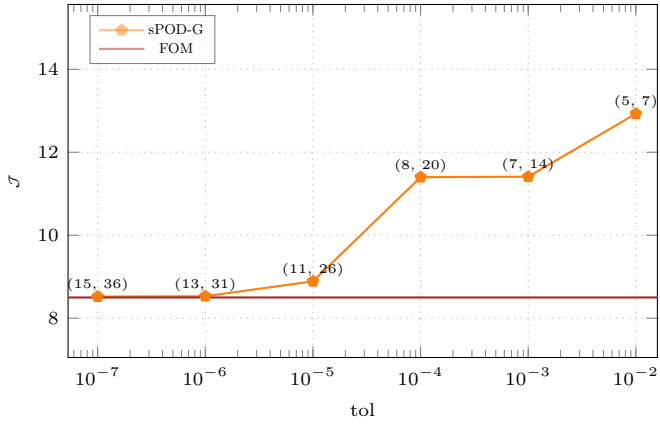
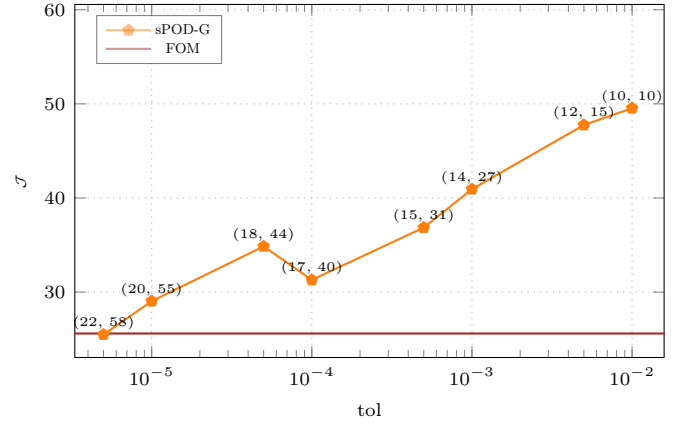
However, we can also use a tolerance-based strategy that prescribes a tolerance which provides a dynamical selection strategy by choosing the number of modes via
\begin{equation}\label{eq_def:epsilon_criteria}
    d = \sum^{\min(n, n_\mathrm{t})}_{i=1} \mathds{1}\left(\frac{\sigma_i}{\sigma_1} > \mathrm{tol}\right)\,
\end{equation}
where $\{\sigma_i\}$ are the singular values for the POD-G or sPOD-G method and the summation counts all values above a relative threshold.
Since the comparison study for POD-G and sPOD-G has already been performed in Figure~\ref{fig:J_shifting}, for the tolerance-based study, we only look at the sPOD-G method and comment on some crucial observations.
In Remark~\ref{rem:control_shape}, we mention that using the control shape functions in a certain way results in an upper bound on the rank of the state snapshot matrix in the stationary frame.
However, this is not guaranteed for the adjoint snapshot matrix in the shifted frame, since the adjoint is the result of the optimality system and cannot be tuned independently.
Thus, it is absolutely possible that the shifted adjoint snapshot matrix is not necessarily a low-rank matrix.
This can be seen from the tolerance-based plot, see Figure~\ref{fig:J_tol_shifting}, where the number of truncated modes for the adjoint is almost always more than the truncated modes for the state for the same tolerance value.
Thus, in principle, one can achieve convergence to the FOM cost functional with fewer modes for the state.

Now, even though the sPOD-G method needs fewer modes than POD-G to reach the cost functional value of the FOM, at this stage, the computational time it takes is not always less than that of the POD-G or the FOM (see Tables~\ref{tab:timing_split_st} and~\ref{tab:timing_split_dt}).
\begin{table}[htp]
  \centering
  \caption{Computational time ($\mathrm{s}$) for crucial steps (Single tilt problem)}
  \label{tab:timing_split_st}
  \begin{tabular}{l|ccc}
    \toprule
    Computational Steps & POD-G ($\mathrm{modes} = 300$) & sPOD-G ($\mathrm{modes} = 25$) & FOM \\
    \midrule
    $n_\mathrm{iter}$  & $5206$  & $1701$  & $1901$ \\
    \midrule
    Basis construction            & $1199.82$  & $752.29$ & $0.00$ \\
    ROM/FOM state solve    & $176.26$  & $124.68$  & $255.69$\\
    Compute $\mathcal{J}$            & $341.13$  & $166.78$ & $76.77$  \\
    ROM/FOM adjoint solve    & $497.02$  & $809.92$ & $369.82$  \\   
    Compute gradient            & $5.60$  & $154.15$ & $11.24$\\
    Update control    & $395.33$  & $451.74$  & $721.55$  \\
    \midrule 
    \textbf{Total}            & $3571.51$  & $2749.00$ & $1729.27$ \\
    \bottomrule
  \end{tabular}
\end{table}
\begin{table}[htp]
  \centering
  \caption{Computational time ($\mathrm{s}$) for crucial steps (Double tilt problem)}
  \label{tab:timing_split_dt}
  \begin{tabular}{l|ccc}
    \toprule
    Computational Steps & POD-G ($\mathrm{modes} = 300$) & sPOD-G ($\mathrm{modes} = 45$) & FOM \\
    \midrule
    $n_\mathrm{iter}$  & $4755$  & $1701$  & $2735$ \\
    \midrule
    Basis construction            & $1106.90$  & $1058.97$ & $0.00$ \\
    ROM/FOM state solve    & $163.86$  & $168.56$  & $363.90$\\
    Compute $\mathcal{J}$            & $310.19$  & $363.07$ & $109.69$  \\
    ROM/FOM adjoint solve    & $449.79$  & $1629.17$ & $521.54$  \\   
    Compute gradient            & $5.27$  & $227.50$ & $14.47$\\
    Update control    & $515.03$  & $972.62$  & $882.80$  \\
    \midrule 
    \textbf{Total}            & $3412.49$  & $4681.64$ & $2304.15$ \\
    \bottomrule
  \end{tabular}
\end{table}
We note here that the FOM is the fastest, since it does not incur the overhead of constructing a reduced basis and also exploits highly optimized sparse-matrix libraries for the state and adjoint equation solution.
Comparing the two ROMs, we see that the sPOD-G method generally takes more time, notably in the reduced-adjoint solve step.
The sPOD-G reduced-adjoint equation is a nonlinear equation, and some of its terms also scale with the FOM dimension in the time integration step. 
This makes it computationally more expensive than the linear POD-G reduced-adjoint equation.
Additionally, we note that, during the basis construction step, we need to assemble shift dependent terms for both the state and adjoint reduced-order equations, which also adds to the total computational time, and this step turns out to be more expensive than that for the POD-G.
However, in our tests, since the POD-G takes more number of optimization steps in total, the cost for the basis construction step is balanced out between the two methods. 
Apart from these, we observe that the sPOD-G reduced-state equation takes more time than the POD-G reduced-state equation in a per-optimization step measure, since the number of modes considered here ($25$ or $45$) are not in principle low-enough to observe any significant speedup.

\section{Conclusions}
For a specific linear hyperbolic optimal control problem, we investigated an FOTR framework for the sPOD method which we compared with its POD analog. It turned out that already for such a relatively simple problem, using the sPOD-G method becomes nontrivial due to the resulting set of nonlinear state and adjoint equations.  
Nevertheless, our numerical results show potential of the approach since the nonlinear reduction approach needs an order of magnitude fewer modes than the POD-G method to capture the FOM dynamics and optimal costs. As a next step, one could study hyperreduction strategies like DEIM \cite{chaturantabut_nonlinear_2010} for the sPOD-G reduced-adjoint equation to speed up its computation. Another possible direction of research could be to dynamically update the shifts in between optimization steps so that the profile could truly be low-rank at every optimization step.

\section*{Acknowledgement}
 We gratefully acknowledge the support of the Deutsche Forschungsgemeinschaft (DFG) as part of the GRK2433 DAEDALUS (DFG Project number: 384950143). 
We also acknowledge the financial support from the SFB TRR154 (DFG project number: 239904186) under sub-project B03.

\bibliographystyle{plain}
\bibliography{abbr,refs}

\end{document}

%% file: State_target_p1.tex
% This file was created with tikzplotlib v0.10.1.
\begin{tikzpicture}

\definecolor{darkgray176}{RGB}{176,176,176}

\begin{groupplot}[group style={group size=2 by 1, horizontal sep=1.2cm}]
\nextgroupplot[
colorbar,
colorbar style={ytick={-0.2,0,0.2,0.4,0.6,0.8,1,1.2},yticklabels={
  \(\displaystyle {\ensuremath{-}0.2}\),
  \(\displaystyle {0.0}\),
  \(\displaystyle {0.2}\),
  \(\displaystyle {0.4}\),
  \(\displaystyle {0.6}\),
  \(\displaystyle {0.8}\),
  \(\displaystyle {1.0}\),
  \(\displaystyle {1.2}\)
},ylabel={}, width=0.1*\pgfkeysvalueof{/pgfplots/parent axis width}, ticklabel style={font=\tiny}},
colormap={mymap}{[1pt]
  rgb(0pt)=(1,1,0.8);
  rgb(1pt)=(1,0.929411764705882,0.627450980392157);
  rgb(2pt)=(0.996078431372549,0.850980392156863,0.462745098039216);
  rgb(3pt)=(0.996078431372549,0.698039215686274,0.298039215686275);
  rgb(4pt)=(0.992156862745098,0.552941176470588,0.235294117647059);
  rgb(5pt)=(0.988235294117647,0.305882352941176,0.164705882352941);
  rgb(6pt)=(0.890196078431372,0.101960784313725,0.109803921568627);
  rgb(7pt)=(0.741176470588235,0,0.149019607843137);
  rgb(8pt)=(0.501960784313725,0,0.149019607843137)
},
height=1.2 * \figureheight,
hide x axis,
hide y axis,
point meta max=1.00000201032067,
point meta min=-1.56173253764891e-05,
tick align=outside,
tick pos=left,
title={\(\displaystyle Q\)},
width=0.7 * \figurewidth,
x grid style={darkgray176},
xmin=0, xmax=800,
xtick style={color=black},
xtick={0,200,400,600,800},
xticklabels={
  \(\displaystyle {0}\),
  \(\displaystyle {200}\),
  \(\displaystyle {400}\),
  \(\displaystyle {600}\),
  \(\displaystyle {800}\)
},
y grid style={darkgray176},
ymin=0, ymax=3360,
ytick style={color=black},
ytick={0,500,1000,1500,2000,2500,3000,3500},
yticklabels={
  \(\displaystyle {0}\),
  \(\displaystyle {500}\),
  \(\displaystyle {1000}\),
  \(\displaystyle {1500}\),
  \(\displaystyle {2000}\),
  \(\displaystyle {2500}\),
  \(\displaystyle {3000}\),
  \(\displaystyle {3500}\)
}
]
\addplot graphics [includegraphics cmd=\pgfimage,xmin=0, xmax=800, ymin=0, ymax=3360] {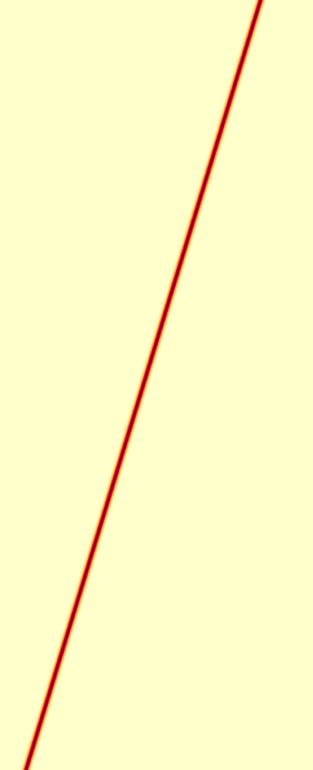};

\nextgroupplot[
colorbar,
colorbar style={ytick={-0.2,0,0.2,0.4,0.6,0.8,1,1.2},yticklabels={
  \(\displaystyle {\ensuremath{-}0.2}\),
  \(\displaystyle {0.0}\),
  \(\displaystyle {0.2}\),
  \(\displaystyle {0.4}\),
  \(\displaystyle {0.6}\),
  \(\displaystyle {0.8}\),
  \(\displaystyle {1.0}\),
  \(\displaystyle {1.2}\)
},ylabel={}, width=0.1*\pgfkeysvalueof{/pgfplots/parent axis width}, ticklabel style={font=\tiny}},
colormap={mymap}{[1pt]
  rgb(0pt)=(1,1,0.8);
  rgb(1pt)=(1,0.929411764705882,0.627450980392157);
  rgb(2pt)=(0.996078431372549,0.850980392156863,0.462745098039216);
  rgb(3pt)=(0.996078431372549,0.698039215686274,0.298039215686275);
  rgb(4pt)=(0.992156862745098,0.552941176470588,0.235294117647059);
  rgb(5pt)=(0.988235294117647,0.305882352941176,0.164705882352941);
  rgb(6pt)=(0.890196078431372,0.101960784313725,0.109803921568627);
  rgb(7pt)=(0.741176470588235,0,0.149019607843137);
  rgb(8pt)=(0.501960784313725,0,0.149019607843137)
},
height=1.2 * \figureheight,
hide x axis,
hide y axis,
point meta max=1.0000020161487,
point meta min=-1.56173253764891e-05,
tick align=outside,
tick pos=left,
title={\(\displaystyle Q_\mathrm{d}\)},
width=0.7 * \figurewidth,
x grid style={darkgray176},
xmin=0, xmax=800,
xtick style={color=black},
xtick={0,200,400,600,800},
xticklabels={
  \(\displaystyle {0}\),
  \(\displaystyle {200}\),
  \(\displaystyle {400}\),
  \(\displaystyle {600}\),
  \(\displaystyle {800}\)
},
y grid style={darkgray176},
ymin=0, ymax=3360,
ytick style={color=black},
ytick={0,500,1000,1500,2000,2500,3000,3500},
yticklabels={
  \(\displaystyle {0}\),
  \(\displaystyle {500}\),
  \(\displaystyle {1000}\),
  \(\displaystyle {1500}\),
  \(\displaystyle {2000}\),
  \(\displaystyle {2500}\),
  \(\displaystyle {3000}\),
  \(\displaystyle {3500}\)
}
]
\addplot graphics [includegraphics cmd=\pgfimage,xmin=0, xmax=800, ymin=0, ymax=3360] {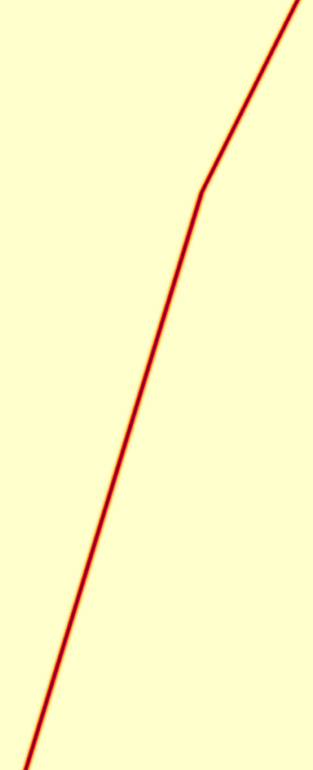};

\end{groupplot}

\draw ({$(current bounding box.south west)!-0.05!(current bounding box.south east)$}|-{$(current bounding box.south west)!0.3!(current bounding box.north west)$}) node[
  scale=0.96,
  anchor=west,
  text=black,
  rotate=90.0
]{time $t$};
\draw ({$(current bounding box.south west)!0.5!(current bounding box.south east)$}|-{$(current bounding box.south west)!-0.2!(current bounding box.north west)$}) node[
  scale=0.96,
  anchor=south,
  text=black,
  rotate=0.0
]{space $x$};
\end{tikzpicture}

%% file: State_target_p4.tex
% This file was created with tikzplotlib v0.10.1.
\begin{tikzpicture}

\definecolor{darkgray176}{RGB}{176,176,176}

\begin{groupplot}[group style={group size=2 by 1, horizontal sep=1.2cm}]
\nextgroupplot[
colorbar,
colorbar style={ytick={-0.2,0,0.2,0.4,0.6,0.8,1,1.2},yticklabels={
  \(\displaystyle {\ensuremath{-}0.2}\),
  \(\displaystyle {0.0}\),
  \(\displaystyle {0.2}\),
  \(\displaystyle {0.4}\),
  \(\displaystyle {0.6}\),
  \(\displaystyle {0.8}\),
  \(\displaystyle {1.0}\),
  \(\displaystyle {1.2}\)
},ylabel={}, width=0.1*\pgfkeysvalueof{/pgfplots/parent axis width}, ticklabel style={font=\tiny}},
colormap={mymap}{[1pt]
  rgb(0pt)=(1,1,0.8);
  rgb(1pt)=(1,0.929411764705882,0.627450980392157);
  rgb(2pt)=(0.996078431372549,0.850980392156863,0.462745098039216);
  rgb(3pt)=(0.996078431372549,0.698039215686274,0.298039215686275);
  rgb(4pt)=(0.992156862745098,0.552941176470588,0.235294117647059);
  rgb(5pt)=(0.988235294117647,0.305882352941176,0.164705882352941);
  rgb(6pt)=(0.890196078431372,0.101960784313725,0.109803921568627);
  rgb(7pt)=(0.741176470588235,0,0.149019607843137);
  rgb(8pt)=(0.501960784313725,0,0.149019607843137)
},
height=1.2 * \figureheight,
hide x axis,
hide y axis,
point meta max=1.00000201032067,
point meta min=-1.56173253764891e-05,
tick align=outside,
tick pos=left,
title={\(\displaystyle Q\)},
width=0.7 * \figurewidth,
x grid style={darkgray176},
xmin=0, xmax=800,
xtick style={color=black},
xtick={0,200,400,600,800},
xticklabels={
  \(\displaystyle {0}\),
  \(\displaystyle {200}\),
  \(\displaystyle {400}\),
  \(\displaystyle {600}\),
  \(\displaystyle {800}\)
},
y grid style={darkgray176},
ymin=0, ymax=3360,
ytick style={color=black},
ytick={0,500,1000,1500,2000,2500,3000,3500},
yticklabels={
  \(\displaystyle {0}\),
  \(\displaystyle {500}\),
  \(\displaystyle {1000}\),
  \(\displaystyle {1500}\),
  \(\displaystyle {2000}\),
  \(\displaystyle {2500}\),
  \(\displaystyle {3000}\),
  \(\displaystyle {3500}\)
}
]
\addplot graphics [includegraphics cmd=\pgfimage,xmin=0, xmax=800, ymin=0, ymax=3360] {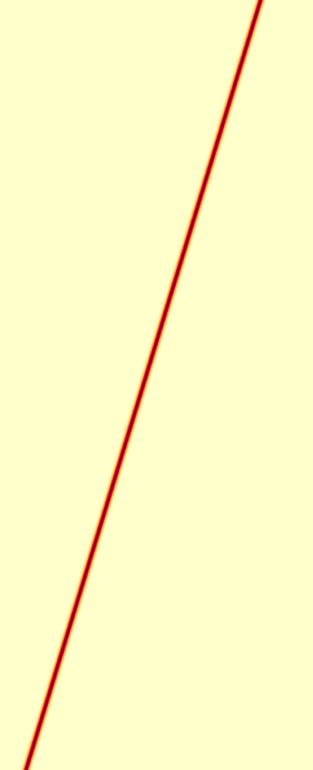};

\nextgroupplot[
colorbar,
colorbar style={ytick={-0.2,0,0.2,0.4,0.6,0.8,1,1.2},yticklabels={
  \(\displaystyle {\ensuremath{-}0.2}\),
  \(\displaystyle {0.0}\),
  \(\displaystyle {0.2}\),
  \(\displaystyle {0.4}\),
  \(\displaystyle {0.6}\),
  \(\displaystyle {0.8}\),
  \(\displaystyle {1.0}\),
  \(\displaystyle {1.2}\)
},ylabel={}, width=0.1*\pgfkeysvalueof{/pgfplots/parent axis width}, ticklabel style={font=\tiny}},
colormap={mymap}{[1pt]
  rgb(0pt)=(1,1,0.8);
  rgb(1pt)=(1,0.929411764705882,0.627450980392157);
  rgb(2pt)=(0.996078431372549,0.850980392156863,0.462745098039216);
  rgb(3pt)=(0.996078431372549,0.698039215686274,0.298039215686275);
  rgb(4pt)=(0.992156862745098,0.552941176470588,0.235294117647059);
  rgb(5pt)=(0.988235294117647,0.305882352941176,0.164705882352941);
  rgb(6pt)=(0.890196078431372,0.101960784313725,0.109803921568627);
  rgb(7pt)=(0.741176470588235,0,0.149019607843137);
  rgb(8pt)=(0.501960784313725,0,0.149019607843137)
},
height=1.2 * \figureheight,
hide x axis,
hide y axis,
point meta max=1.0000020161487,
point meta min=-1.56173253764891e-05,
tick align=outside,
tick pos=left,
title={\(\displaystyle Q_\mathrm{d}\)},
width=0.7 * \figurewidth,
x grid style={darkgray176},
xmin=0, xmax=800,
xtick style={color=black},
xtick={0,200,400,600,800},
xticklabels={
  \(\displaystyle {0}\),
  \(\displaystyle {200}\),
  \(\displaystyle {400}\),
  \(\displaystyle {600}\),
  \(\displaystyle {800}\)
},
y grid style={darkgray176},
ymin=0, ymax=3360,
ytick style={color=black},
ytick={0,500,1000,1500,2000,2500,3000,3500},
yticklabels={
  \(\displaystyle {0}\),
  \(\displaystyle {500}\),
  \(\displaystyle {1000}\),
  \(\displaystyle {1500}\),
  \(\displaystyle {2000}\),
  \(\displaystyle {2500}\),
  \(\displaystyle {3000}\),
  \(\displaystyle {3500}\)
}
]
\addplot graphics [includegraphics cmd=\pgfimage,xmin=0, xmax=800, ymin=0, ymax=3360] {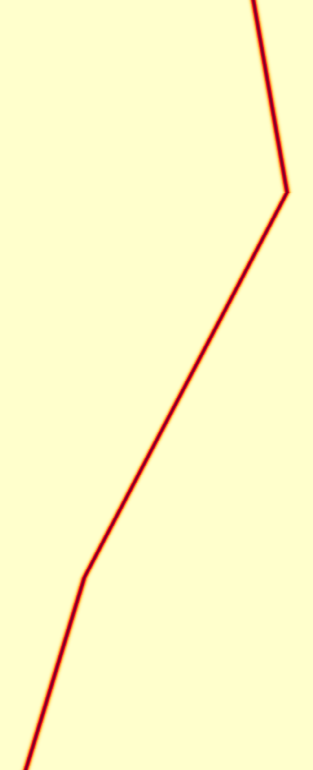};

\end{groupplot}

\draw ({$(current bounding box.south west)!-0.05!(current bounding box.south east)$}|-{$(current bounding box.south west)!0.3!(current bounding box.north west)$}) node[
  scale=0.96,
  anchor=west,
  text=black,
  rotate=90.0
]{time $t$};
\draw ({$(current bounding box.south west)!0.5!(current bounding box.south east)$}|-{$(current bounding box.south west)!-0.2!(current bounding box.north west)$}) node[
  scale=0.96,
  anchor=south,
  text=black,
  rotate=0.0
]{space $x$};
\end{tikzpicture}

%% file: FOTR_modesVsCost_Adaptive_P1.tex
\begin{tikzpicture}

  % Define colors (same as your original)
  \definecolor{brown}{RGB}{165,42,42}
  \definecolor{darkgray176}{RGB}{176,176,176}
  \definecolor{green}{RGB}{0,128,0}
  \definecolor{lightgray204}{RGB}{204,204,204}
  \definecolor{darkorange25512714}{RGB}{255,127,14}
  \definecolor{forestgreen4416044}{RGB}{44,160,44}
  \definecolor{sienna}{RGB}{160,82,45}
  \definecolor{steelblue31119180}{RGB}{31,119,180}
    \definecolor{crimson2143940}{RGB}{214,39,40}

  \begin{axis}[
    height=0.9*\figureheight,
    width=1.2*\figurewidth,
    % log basis x={10},
    % log basis y={10},
    % xmode=log,
    % ymode=log,
    xmin=5.01187233627271e-08,   xmax=305,
    ymin=5.95019163314602,      ymax=100.452607734195,
    xlabel={modes},
    ylabel={$\scriptstyle{\mathcal{J}}$},
    xmajorgrids,
    ymajorgrids,
    x grid style={darkgray176},
    y grid style={darkgray176},
    xtick={10,30,50,100,150,200,250,300,350,400,450,500},
    xticklabels={
      \(\displaystyle 10\),
      \(\displaystyle 30\),
      \(\displaystyle 50\),
      \(\displaystyle 100\),
      \(\displaystyle 150\),
      \(\displaystyle 200\),
      \(\displaystyle 250\),
      \(\displaystyle 300\),
      \(\displaystyle 350\),
      \(\displaystyle 400\),
      \(\displaystyle 450\),
      \(\displaystyle 500\),
    },
    ytick={0,10,20,30,40,50,60,70,80,90},
    yticklabels={
      \(\displaystyle 0\),
      \(\displaystyle 10\),
      \(\displaystyle 20\),
      \(\displaystyle 30\),
      \(\displaystyle 40\),
      \(\displaystyle 50\),
      \(\displaystyle 60\),
      \(\displaystyle 70\),
      \(\displaystyle 80\),
      \(\displaystyle 90\),
    },
    legend style={
    font=\footnotesize,
      fill opacity=0.5,
      draw opacity=0.5,
      text opacity=1,
      nodes={scale=0.6},
      draw=lightgray204,
      at={(0.95,0.95)}
    }
  ]

    %––––––––––––––––––––––––––––––––––––––––––––––––––––––––––––––
    % Series 1: FOTR primal basis, with minimal‐size labels
    \addplot+[
      thick,
      steelblue31119180,
      mark=ball,
      mark size=2,
      mark options={solid},
      nodes near coords,
      every node near coord/.append style={
        font=\tiny,
        inner sep=0pt,
        outer sep=0pt,
        yshift=2pt
      },
      point meta=explicit symbolic
    ]
    table[
      row sep=\\,
      x=x,
      y=y,
      meta=label
    ] {
      x           y                    label \\
5   38.15405614976782    {} \\
10   38.15405614976782    {} \\
20   38.15405614976782    {} \\
30   93.7378234605576    {} \\
40   27.76386732586036    {} \\
50   23.643730960202465   {} \\
60   23.428035654877355    {} \\
70   24.146577866222717    {} \\
80   18.609110714082213    {} \\
90   21.146653565464113    {} \\
100   18.219358322239636    {} \\
200   22.468443108293613    {} \\
300   8.499739040747947    {} \\
400   8.50106107669594    {} \\
500   8.529381822095711    {} \\
    };
    \addlegendentry{\small POD-G}
    
    %––––––––––––––––––––––––––––––––––––––––––––––––––––––––––––––
    % Series 2: FOTR primal basis, with minimal‐size labels
    \addplot+[
      thick,
      darkorange25512714,
      mark=pentagon*,
      mark size=2,
      mark options={solid},
      nodes near coords,
      every node near coord/.append style={
        font=\tiny,
        inner sep=0pt,
        outer sep=0pt,
        yshift=2pt
      },
      point meta=explicit symbolic
    ]
    table[
      row sep=\\,
      x=x,
      y=y,
      meta=label
    ] {
      x           y                    label \\
2   20.932194549908225    {} \\
5   14.984190497196803    {} \\
8   12.620877355901579    {} \\
10   11.32732221477353    {} \\
12   8.920298948635645    {} \\
15   8.573004652405182    {} \\
20   8.633966170133734    {} \\
25   8.523283158454342    {} \\
30   8.549935559067082    {} \\
    };
    \addlegendentry{\small sPOD-G}

    %––––––––––––––––––––––––––––––––––––––––––––––––––––––––––––––
    % Series 3: FOM horizontal line (no labels needed here)
    \addplot[
    thick,
      mark=none,
      brown,
      samples=150
    ] coordinates {
      (5.01187233627271e-08,8.499)
      (510.0199526231496888,8.499)
    };
    \addlegendentry{\small FOM}

  \end{axis}
\end{tikzpicture}

%% file: FOTR_modesVsCost_Adaptive_P.tex
\begin{tikzpicture}

  % Define colors (same as your original)
  \definecolor{brown}{RGB}{165,42,42}
  \definecolor{darkgray176}{RGB}{176,176,176}
  \definecolor{green}{RGB}{0,128,0}
  \definecolor{lightgray204}{RGB}{204,204,204}
  \definecolor{darkorange25512714}{RGB}{255,127,14}
  \definecolor{forestgreen4416044}{RGB}{44,160,44}
  \definecolor{sienna}{RGB}{160,82,45}
  \definecolor{steelblue31119180}{RGB}{31,119,180}
    \definecolor{crimson2143940}{RGB}{214,39,40}

  \begin{axis}[
    height=0.9*\figureheight,
    width=1.2*\figurewidth,
    % log basis x={10},
    % log basis y={10},
    % % xmode=log,
    % ymode=log,
    xmin=5.01187233627271e-08,   xmax=305,
    ymin=15.95019163314602,      ymax=142.452607734195,
    xlabel={modes},
    ylabel={$\scriptstyle{\mathcal{J}}$},
    xmajorgrids,
    ymajorgrids,
    x grid style={darkgray176},
    y grid style={darkgray176},
    xtick={10,30,50,100,150,200,250,300,350,400,450,500},
    xticklabels={
      \(\displaystyle 10\),
      \(\displaystyle 30\),
      \(\displaystyle 50\),
      \(\displaystyle 100\),
      \(\displaystyle 150\),
      \(\displaystyle 200\),
      \(\displaystyle 250\),
      \(\displaystyle 300\),
      \(\displaystyle 350\),
      \(\displaystyle 400\),
      \(\displaystyle 450\),
      \(\displaystyle 500\),
    },
    ytick={0,10,20,30,40,50,60,70,80,90,100,110,120,130,140},
    yticklabels={
      \(\displaystyle 0\),
      \(\displaystyle 10\),
      \(\displaystyle 20\),
      \(\displaystyle 30\),
      \(\displaystyle 40\),
      \(\displaystyle 50\),
      \(\displaystyle 60\),
      \(\displaystyle 70\),
      \(\displaystyle 80\),
      \(\displaystyle 90\),
      \(\displaystyle 100\),
      \(\displaystyle 110\),
      \(\displaystyle 120\),
      \(\displaystyle 130\),
      \(\displaystyle 140\),
    },
    legend style={
    font=\footnotesize,
      fill opacity=0.5,
      draw opacity=0.5,
      text opacity=1,
      nodes={scale=0.6},
      draw=lightgray204,
      at={(0.95,0.95)}
    }
  ]

    %––––––––––––––––––––––––––––––––––––––––––––––––––––––––––––––
    % Series 1: FOTR primal basis, with minimal‐size labels
    \addplot+[
      thick,
      steelblue31119180,
      mark=ball,
      mark size=2,
      mark options={solid},
      nodes near coords,
      every node near coord/.append style={
        font=\tiny,
        inner sep=0pt,
        outer sep=0pt,
        yshift=2pt
      },
      point meta=explicit symbolic
    ]
    table[
      row sep=\\,
      x=x,
      y=y,
      meta=label
    ] {
      x           y                    label \\
5   120.56404171231596    {} \\
10   120.56404171231596    {} \\
20   108.40299865772297    {} \\
30   87.32349740651303    {} \\
40   74.63183444235707    {} \\
50   71.53342894728472   {} \\
60   72.3604986025728    {} \\
70   71.57590726695693    {} \\
80   71.6115485926935    {} \\
90   60.81783652543514    {} \\
100   67.47097239594574    {} \\
200   70.49582154552034    {} \\
300   25.41968411559036    {} \\
400   25.418627831960162    {} \\
500   25.41861049725641    {} \\
    };
    \addlegendentry{\small POD-G}
    
    %––––––––––––––––––––––––––––––––––––––––––––––––––––––––––––––
    % Series 2: FOTR primal basis, with minimal‐size labels
    \addplot+[
      thick,
      darkorange25512714,
      mark=pentagon*,
      mark size=2,
      mark options={solid},
      nodes near coords,
      every node near coord/.append style={
        font=\tiny,
        inner sep=0pt,
        outer sep=0pt,
        yshift=2pt
      },
      point meta=explicit symbolic
    ]
    table[
      row sep=\\,
      x=x,
      y=y,
      meta=label
    ] {
      x           y                    label \\
2   71.59203212823007    {} \\
5   63.997089415523455    {} \\
8   50.307273129094476    {} \\
10   46.57595366736163    {} \\
12   46.29794142998056    {} \\
15   48.008380350750905    {} \\
20   29.277032827720873    {} \\
25   26.903461024025088    {} \\
30   34.96886710571803    {} \\
35   28.04925344566443    {} \\
40   28.16920937938442    {} \\
45   25.6809988578256    {} \\
50   25.511550225684015    {} \\
    };
    \addlegendentry{\small sPOD-G}

    %––––––––––––––––––––––––––––––––––––––––––––––––––––––––––––––
    % Series 3: FOM horizontal line (no labels needed here)
    \addplot[
    thick,
      mark=none,
      brown,
      samples=150
    ] coordinates {
      (5.01187233627271e-08,25.60)
      (510.0199526231496888,25.60)
    };
    \addlegendentry{\small FOM}

  \end{axis}
\end{tikzpicture}

%% file: FOTR_tolVsCost_Adaptive_P1.tex
\begin{tikzpicture}

  % Define colors (same as your original)
  \definecolor{brown}{RGB}{165,42,42}
  \definecolor{darkgray176}{RGB}{176,176,176}
  \definecolor{green}{RGB}{0,128,0}
  \definecolor{lightgray204}{RGB}{204,204,204}

  \definecolor{crimson2143940}{RGB}{214,39,40}
  \definecolor{darkorange25512714}{RGB}{255,127,14}
  \definecolor{forestgreen4416044}{RGB}{44,160,44}
  \definecolor{sienna}{RGB}{160,82,45}
  \definecolor{steelblue31119180}{RGB}{31,119,180}

  \begin{axis}[
    height=0.9*\figureheight,
    width=1.2*\figurewidth,
    % log basis x={10},
    % log basis y={10},
    xmode=log,
    % ymode=log,
    xmin=5.30957344480193e-8, xmax=0.0158489319246111,
    ymin=7.052039401636526,    ymax=15.56193659558,
    xlabel={tol},
    ylabel={$\scriptstyle{\mathcal{J}}$},
    xmajorgrids,
    ymajorgrids,
    x grid style={darkgray176},
    y grid style={darkgray176},
    xtick={1e-10,1e-09,1e-08,1e-07,1e-06,1e-05,1e-04,1e-03,1e-02,1e-01,1},
    xticklabels={
      \(\displaystyle 10^{-10}\),
      \(\displaystyle 10^{-9}\),
      \(\displaystyle 10^{-8}\),
      \(\displaystyle 10^{-7}\),
      \(\displaystyle 10^{-6}\),
      \(\displaystyle 10^{-5}\),
      \(\displaystyle 10^{-4}\),
      \(\displaystyle 10^{-3}\),
      \(\displaystyle 10^{-2}\),
      \(\displaystyle 10^{-1}\),
      \(\displaystyle 10^{0}\)
    },
    ytick={0,2,4,6,8,10,12,14,16,18,20,22,24},
    yticklabels={
      \(\displaystyle 0\),
      \(\displaystyle 2\),
      \(\displaystyle 4\),
      \(\displaystyle 6\),
      \(\displaystyle 8\),
      \(\displaystyle 10\),
      \(\displaystyle 12\),
      \(\displaystyle 14\),
      \(\displaystyle 16\),
      \(\displaystyle 18\),
      \(\displaystyle 20\),
      \(\displaystyle 22\),
      \(\displaystyle 24\),
    },
    legend style={
    font=\footnotesize,
      fill opacity=0.5,
      draw opacity=0.5,
      text opacity=1,
      nodes={scale=0.6},
      at={(0.25,0.97)},
      anchor=north east
    }
  ]

    %––––––––––––––––––––––––––––––––––––––––––––––––––––––––––––––
    % Series 1: FOTR primal basis, with “label” column
    \addplot+[
      thick,
      darkorange25512714,
      mark=pentagon*,
      mark size=2,
      mark options={solid},
      nodes near coords,
      every node near coord/.append style={
        font=\tiny,
        inner sep=0pt,
        outer sep=0pt,
        yshift=2pt
      },
      point meta=explicit symbolic
    ]
    table[
      row sep=\\,
      x=x,
      y=y,
      meta=label
    ] {
      x         y                   label \\
    % 1e-10   8.517447201467098    {} \\
    % 5e-10   8.526703189296573    {} \\
    % 1e-09   8.610388351752109    {} \\
    % 5e-09   8.612578034843901    {} \\
    % 1e-08   8.541989893157576      {} \\
    % 5e-08  8.828110559703173    {} \\
    1e-07  8.519439899757808      {\textcolor{black}{(15, 36)}} \\
    % 5e-07   8.518832760410726      {} \\
    1e-06   8.528361811754714    {\textcolor{black}{(13, 31)}} \\
    % 5e-06   8.942737320709917    {} \\
    1e-05   8.887788115224838    {\textcolor{black}{(11, 26)}} \\
    % 5e-05   9.151878880266441      {} \\
    0.0001  11.399311063899997    {\textcolor{black}{(8, 20)}} \\
    % 0.0005  10.78779457319443      {} \\
    0.001   11.410118711178116      {\textcolor{black}{(7, 14)}} \\
    % 0.005   14.401316307676952     {} \\
    0.01    12.927934212155177     {\textcolor{black}{(5, 7)}} \\
    };
    \addlegendentry{\small sPOD-G}

    %––––––––––––––––––––––––––––––––––––––––––––––––––––––––––––––
    % Series 3: FOM horizontal line (no labels needed)
    \addplot[
    thick,
      mark=none,
      brown,
      samples=150
    ] coordinates {
      (5.01187233627271e-11, 8.499)
      (0.0199526231496888,   8.499)
    };
    \addlegendentry{\small FOM}

  \end{axis}
\end{tikzpicture}

%% file: FOTR_tolVsCost_Adaptive_P.tex
\begin{tikzpicture}

  % Define colors (same as your original)
  \definecolor{brown}{RGB}{165,42,42}
  \definecolor{darkgray176}{RGB}{176,176,176}
  \definecolor{green}{RGB}{0,128,0}
  \definecolor{lightgray204}{RGB}{204,204,204}

  \definecolor{crimson2143940}{RGB}{214,39,40}
  \definecolor{darkorange25512714}{RGB}{255,127,14}
  \definecolor{forestgreen4416044}{RGB}{44,160,44}
  \definecolor{sienna}{RGB}{160,82,45}
  \definecolor{steelblue31119180}{RGB}{31,119,180}

  \begin{axis}[
    height=0.9*\figureheight,
    width=1.2*\figurewidth,
    % log basis x={10},
    % log basis y={10},
    xmode=log,
    % ymode=log,
    xmin=3.30957344480193e-06, xmax=0.0158489319246111,
    ymin=23.052039401636526,    ymax=60.56193659558,
    xlabel={tol},
    ylabel={$\scriptstyle{\mathcal{J}}$},
    xmajorgrids,
    ymajorgrids,
    x grid style={darkgray176},
    y grid style={darkgray176},
    xtick={1e-10,1e-09,1e-08,1e-07,1e-06,1e-05,1e-04,1e-03,1e-02,1e-01,1},
    xticklabels={
      \(\displaystyle 10^{-10}\),
      \(\displaystyle 10^{-9}\),
      \(\displaystyle 10^{-8}\),
      \(\displaystyle 10^{-7}\),
      \(\displaystyle 10^{-6}\),
      \(\displaystyle 10^{-5}\),
      \(\displaystyle 10^{-4}\),
      \(\displaystyle 10^{-3}\),
      \(\displaystyle 10^{-2}\),
      \(\displaystyle 10^{-1}\),
      \(\displaystyle 10^{0}\)
    },
    ytick={0,10,20,30,40,50,60,70,80,90,100,110,120},
    yticklabels={
      \(\displaystyle 0\),
      \(\displaystyle 10\),
      \(\displaystyle 20\),
      \(\displaystyle 30\),
      \(\displaystyle 40\),
      \(\displaystyle 50\),
      \(\displaystyle 60\),
      \(\displaystyle 70\),
      \(\displaystyle 80\),
      \(\displaystyle 90\),
      \(\displaystyle 100\),
      \(\displaystyle 110\),
      \(\displaystyle 120\),
    },
    legend style={
    font=\footnotesize,
      fill opacity=0.5,
      draw opacity=0.5,
      text opacity=1,
      nodes={scale=0.6},
      at={(0.25,0.97)},
      anchor=north east
    }
  ]

    %––––––––––––––––––––––––––––––––––––––––––––––––––––––––––––––
    % Series 1: FRTO primal basis, with “label” column
    \addplot+[
      thick,
      darkorange25512714,
      mark=pentagon*,
      mark size=2,
      mark options={solid},
      nodes near coords,
      every node near coord/.append style={
        font=\tiny,
        inner sep=0pt,
        outer sep=0pt,
        yshift=2pt
      },
      point meta=explicit symbolic
    ]
    table[
      row sep=\\,
      x=x,
      y=y,
      meta=label
    ] {
      x         y                   label \\
    % 1e-10   44.08042967173336    {} \\
    % 5e-10   25.481747160142334    {} \\
    % 1e-09   25.48172346229962    {} \\
    % 5e-09   25.4817218344864    {} \\
    % 1e-08   25.482645417181665      {} \\
    % 5e-08  25.531680240060382    {} \\
    % 1e-07  25.48370341847824      {} \\
    % 5e-07   25.52591317597811      {} \\
    % 1e-06   25.95864136096509    {} \\
    5e-06   25.486161204432893    {\textcolor{black}{(22, 58)}} \\
    1e-05   29.036737892916122    {\textcolor{black}{(20, 55)}} \\
    5e-05   34.852199564291176      {\textcolor{black}{(18, 44)}} \\
    0.0001  31.293636643599704    {\textcolor{black}{(17, 40)}} \\
    0.0005  36.85510104040199      {\textcolor{black}{(15, 31)}} \\
    0.001   40.93272663257237      {\textcolor{black}{(14, 27)}} \\
    0.005   47.75483635105342     {\textcolor{black}{(12, 15)}} \\
    0.01    49.538307942658804     {\textcolor{black}{(10, 10)}} \\
    };
    \addlegendentry{\small sPOD-G}

    %––––––––––––––––––––––––––––––––––––––––––––––––––––––––––––––
    % Series 2: FOM horizontal line (no labels needed)
    \addplot[
    thick,
      mark=none,
      brown,
      samples=150
    ] coordinates {
      (5.01187233627271e-11, 25.60)
      (0.0199526231496888,   25.60)
    };
    \addlegendentry{\small FOM}

  \end{axis}
\end{tikzpicture}